\documentclass[
twocolumn,
tightenlines,
10pt,
longbibliography,
showpacs,
nofootinbib,
notitlepage,
superscriptaddress]{revtex4-2}

\usepackage{graphicx}
\usepackage{dcolumn}
\usepackage{bm}
\usepackage{lipsum}  
\usepackage{xcolor}
\usepackage{soul}
\usepackage{amsmath, amssymb, amsthm, braket, bbold, enumitem, dsfont, relsize, physics}
\usepackage{makecell, graphicx, subfiles, multirow, wrapfig, tabularx}
\usepackage[font=scriptsize]{subcaption}
\usepackage[font=scriptsize, justification=raggedright]{caption}
\usepackage[page]{appendix}
\usepackage[colorlinks, allcolors=black]{hyperref}
\usepackage[capitalize, nameinlink]{cleveref}
\usepackage{tikz}
\usetikzlibrary{quantikz}
\usetikzlibrary{shapes.geometric, arrows}
\usetikzlibrary{positioning}
\setlist[itemize]{leftmargin=10pt}

\newcommand{\mi}{\mathrm{i}} 

\begin{document}

\title{Lagrangian Duality in Quantum Optimization: \\Overcoming QUBO Limitations for Constrained Problems}

\author{Einar Gabbassov}
\affiliation{1QB Information Technologies (1QBit), Vancouver, BC, Canada}
\affiliation{Department of Applied Mathematics, University of Waterloo, Waterloo, ON, Canada}

\author{Gili Rosenberg}
\affiliation{1QB Information Technologies (1QBit), Vancouver, BC, Canada}
\affiliation{Amazon Quantum Solutions Lab, Seattle, WA 98170, USA}
\thanks{This work was done prior to joining Amazon.}

\author{Artur Scherer}
\affiliation{1QB Information Technologies (1QBit), Vancouver, BC, Canada}

\date{April 22, 2024}
\begin{abstract}
    We propose an approach to solving constrained combinatorial optimization problems based on embedding the concept of Lagrangian duality into the framework of adiabatic quantum computation. Within the setting of circuit-model fault-tolerant quantum computation, we demonstrate that this approach achieves a quadratic improvement in circuit depth and maintains a constraint-independent circuit width in contrast to the prevalent approach of solving constrained problems via reformulations based on the quadratic unconstrained
    binary optimization (QUBO) framework. Our study includes a detailed review of the limitations encountered when using QUBO for constrained optimization. We show that the proposed method overcomes these limitations by encoding the optimal solution at an energetically elevated level of a simpler problem Hamiltonian, which results in substantially more resource-efficient quantum circuits. We consolidate our strategy  with a detailed analysis on how the concepts of Lagrangian duality such as duality gap and complementary slackness relate to the success probability of sampling the optimal solution. Our findings are illustrated by benchmarking the Lagrangian dual approach against the QUBO approach using the NP-complete binary knapsack problem.
\end{abstract}

\maketitle

\section{\label{sec:introduction}Introduction}
Solving optimization problems is considered to be one of the most important potential 
applications of quantum computers. A variety of approaches have been proposed, including quantum algorithms requiring a fault-tolerant quantum computer for their implementation, and quantum heuristics such as quantum annealing~\cite{morita2008mathematical} and QAOA~\cite{farhi2014quantum}. 
Real-world optimization problems frequently also include constraints. Normally, quantum annealing and QAOA only apply to unconstrained problems, which limits their practical applicability; yet it is possible to recast constrained problems as  unconstrained ones to which then the established methods can be applied. The most common approach is based on reformulating a linear-constrained problem into a QUBO problem \cite{lewis2017quadratic}. Essentially, QUBO incorporates linear constraints as quadratic penalties into an objective function. However, such reformulations give rise to other challenges. For example, quadratic penalties for inequality constraints require a large number of additional slack variables; these slack variables significantly increase the search space of a problem and make the optimization landscape more rugged. Moreover, converting constraints into squared penalty terms often results in all-to-all connectivity between a problem's variables, making a problem more complex. Therefore, the QUBO approach to solving constrained problems 
significantly increases the complexity of the underlying algorithm entailing a substantial overhead of quantum resources.

To address the challenges outlined above, we introduce a more resource-efficient strategy which combines the theory of Lagrangian duality \cite{geoffrion1974lagrangean} and the concept of adiabatic quantum computation (AQC) \cite{farhi2014quantum}. The proposed approach places the optimal solution of the original problem within a bounded, energetically elevated eigensubspace of a specially designed Hamiltonian. This contrasts with traditional AQC methods that encode the optimal solution in the lowest-energy eigenstate. 
To reach this energetically elevated subspace effectively, it is necessary to relax the conditions of adiabaticity \cite{farhi2014quantum}. This relaxation implies that the evolution time need not be excessively prolonged, allowing for a more flexible and potentially faster approach to reaching the desired state.

The energetically elevated eigensubspace contains the optimal solution of the original problem  whose separation from the ground state energy is upper-bounded by a predetermined quantity $\epsilon \geq 0$ that bounds the Lagrangian duality gap. We demonstrate that evolving the system from the ground state of an initial Hamiltonian effectively prepares a state that has a significant overlap with the eigenstate encoding the optimal solution. Furthermore, we show that the extent of this overlap is determined by $\epsilon$, and that for any $\epsilon$, there exist an initial Hamiltonian and a total evolution time $T$ that maximize the overlap.

Capitalizing on these findings, we formulate a more efficient approach to tackling constrained combinatorial optimization (CO) problems. Our approach is based on Discretized Adiabatic Quantum Computation (DAQC) tailored to exploiting  Lagrangian duality principles. We refer to this framework as "Lagrangian dual DAQC" (LD-DAQC).

LD-DAQC incorporates a Trotterized approximation of adiabatic evolution, thus facilitating a discrete implementation of AQC (hence the name DAQC). This results in quantum circuit models that are reminiscent of the renowned Quantum Approximate Optimization Algorithm (QAOA). However, a key distinction lies in DAQC's independence from the classical optimization feedback loop required for optimal circuit parameterization in QAOA. Instead, DAQC leverages an approximation of local adiabatic evolution, a strategy previously established to yield quadratic speedup in search problems \cite{roland2002quantum}. This approximation is realized through a pre-defined adiabatic schedule.

Despite these differences, our LD-DAQC approach remains compatible with QAOA and other circuit-model quantum algorithms focused on variational quantum state preparation. Moreover, this approach is equally applicable to analog AQC and quantum annealing scenarios, where the discretization step is not necessary.

The contributions of this study are as follows:
\begin{itemize}
   \item We propose an adiabatic quantum computation  approach to solving constrained optimization problems  based on Lagrangian duality and show that this approach is simpler and more efficient than QUBO-based approaches. Furthermore, we establish a theoretical relation between Lagrangian duality and the probability of sampling an optimal solution to an optimization problem.
   \item Using a circuit-model framework realized by DAQC, we demonstrate a quadratic improvement in circuit complexity and runtime over a QUBO-based approach for problems with linear inequality constraints.
   \item We show that LD-DAQC circuits have significantly lower connectivity requirements than QUBO-DAQC or QUBO-based QAOA circuits. Additionally, the new approach eliminates the necessity for ancillary qubits, negating the need for logarithmic overheads in qubit count arising from constraint slack variables.
   
   \item For a concrete analysis of QUBO-based and Lagrangian-based formulations and their corresponding quantum circuits, we use the NP-complete binary knapsack problem (KP) as a reference problem for comparison.
   \item Our numerical study demonstrates that 
   LD-DAQC outperforms QUBO-DAQC on a large dataset of KPs of different problem sizes. Crucially, we demonstrate the consistent performance of LD-DAQC irrespective of the values in a problem's coefficients specifying the constraints. Conversely, QUBO-based methods display a polylogarithmic performance degradation as a problem's coefficients become larger.
\end{itemize}
This paper is structured as follows. In~\cref{sec:qubo_issues_main}, we first review the inherent challenges associated with the QUBO approach to solving constrained problems, from both classical optimization and quantum algorithms perspectives. \cref{sec:lagrangian_dual_approach} introduces the classical theory of Lagrangian duality and shows how this theory is integrated within the farmework of adiabatic quantum computation. \cref{sec:circuit} presents a detailed construction of an efficient LD-DAQC circuit along with  an analysis of its complexity. Lastly, \cref{sec:numerical} and ~\cref{sec:exp_setup} illustrate the benefits of the proposed methodology through a series of numerical experiments.

\subsection{\label{sec:prev_work}Previous Work}
Our optimization scheme is based on adiabatic quantum computation (AQC)~\cite{albash2018adiabatic} which is polynomially equivalent to the standard circuit-model of quantum computation~\cite{aharonov2008adiabatic}.
A number of meta-heuristics inspired by AQC have been extensively explored for solving combinatorial optimization problems on quantum devices~\cite{gabbassov2022transit, morita2008mathematical, sankar2021benchmark, gilyen2019optimizing, apolloni1989quantum}. The most notable is quantum annealing~\cite{morita2008mathematical} which exploits decreasing quantum fluctuations to search for a lowest-energy eigenstate (ground state) of a Hamiltonian that encodes a given combinatorial optimization problem. For example, in the spin-glass annealing paradigm, this heuristic uses a time-dependent Hamiltonian consisting of two non-commuting terms: a transverse field Hamiltonian that is used for the mixing process, and an Ising Hamiltonian whose ground state encodes the solution to the problem.  The transverse field Hamiltonian is gradually attenuated, while the Ising Hamiltonian is gradually amplified during the annealing process. Numerous subsequent works, e.g.~\cite{seki2012quantum, galindo2020optimal}, focused on determining suitable scheduling strategies for mixing to find the ground state with higher probability.

Variational approaches for determining AQC scheduling are discussed in \cite{matsuura2020vanqver, matsuura2021variationally}. Generally, variational methods involve an outer classical feedback loop that uses  continuous optimization techniques to find a suitable parametrization of a schedule. QAOA, which has commonly become associated with a NISQ-type algorithm, and which can be viewed as a diabatic counterpart to discretized AQC, is also based on a variational hybrid \mbox{quantum--classical} protocol to optimize the gate parameters. 
Although such variational protocols allow for a greater flexibility, they come with severe challenges, see, e.g., \cite{sankar2021benchmark}. In particular, they require a large number of runs of the scheduled quantum evolution on a quantum device. Moreover, most variational approaches optimize for the expected energy of a system rather than the probability of sampling an optimal solution \cite{mcclean2016theory}. If the energy function is not convex with respect to the variational parameters, then the task of finding optimal parameters is in itself an NP-hard problem. Moreover, even if a suitable schedule is determined for some problem instance, each new problem instance requires finding a new parametrization.

The use of Lagrangian duality within the context of quantum annealing has been considered in~\cite{ronagh2016solving, yonaga2020solving,}. The study in \cite{ronagh2016solving} presents a method for solving the Lagrangian dual of a binary quadratic programming problem with inequality constraints. The proposed method successfully integrates the Lagrangian duality with branch-and-bound and quantum annealing heuristics. The work in~\cite{karimi2017subgradient} considers a quantum subgradient method for finding an optimal primal-dual pair for the Lagrangian dual of a constrained binary problem. The subgradients are computed using a quantum annealer and then used in a classical descent algorithm. However, to the best of our knowledge, the use of Lagrangian duality in the context of circuit-model quantum computation has not been previously explored. In particular, its implementation by a quantum circuit and the associated efficiency have not been examined in previous studies.

Both quantum annealing and QAOA are well-suited for optimizing unconstrained binary problems such as \mbox{\textsc{MaxCut}} which belongs to the class of QUBO problems. Constrained problems are canonically reformulated as QUBO problems by integrating the constraints into an objective function as quadratic penalties. The complications introduced by quadratic penalties are addressed in~\cite{hen2016quantum}; to avoid them, it is suggested to utilize a suitably tailored driver Hamiltonian that commutes with the operator representing the linear equality constraint and initialize the system in the distinct state that is simultaneously a ground state of the initial Hamiltonian and an eigenstate of the constraint term. While this approach allows one to omit quadratic penalties for equality constraints, it becomes increasingly hard to prepare a suitable initial state in the presence of multiple equality constraints. This approach can be generalized to inequality constraints by introducing slack variables. However, as shown in~\cref{sec:qubo_issues_main}, the number of slack variables has an explicit dependence on the constraint bound and implicit dependence on the problem size. This renders all approaches that involve slack variables very costly, as it not only implies a large overhead of quantum resources but also increases the search space of a problem.

\section{\label{sec:qubo_issues_main} Challenges of the QUBO approach}
In the context of quantum optimization, the most common approach to solving a constrained combinatorial problem is to reformulate the constrained problem into a QUBO problem and then recast the resulting QUBO problem into the problem of finding the ground state of the corresponding problem Hamiltonian. The QUBO model has been successful and popular in quantum and quantum-inspired optimization. It is the go-to model for many classical and quantum optimization hardware producers. Indeed, most of the cutting-edge optimization hardware platforms implement QUBO solvers. For example, devices such as D-Wave’s quantum annealer~\cite{boixo2014evidence}, NTT’s coherent Ising machine~\cite{inagaki2016coherent}, Fujitsu’s digital annealer~\cite{aramon2019physics} and Toshiba’s simulated quantum bifurcation machine~\cite{goto2019combinatorial}, are all designed to solve QUBO problems. This popularity is due to the equivalence between QUBO and the Ising model \cite{brush1967history}. Platform-agnostic quantum algorithms like AQC and QAOA also require a problem to be reformulated as QUBO or as a more general PUBO which refers to polynomial unconstrained binary optimization \cite{glover2011polynomial}.

While QUBO is a prevalent model used in quantum or quantum-inspired optimization, it has severe shortcomings when applied to constrained problems. Reformulation of a constrained problem into a QUBO problem 
often results in a substantial resource overhead and requirements of 
high qubit connectivity, which in the framework of circuit-model quantum computation implies increased circuit width and depth. For example, all inequality constraints must be incorporated into an objective function as quadratic penalties. Such quadratic penalties 
usually introduce an ill-behaved optimization landscape and 
may require many additional auxiliary binary variables 
that significantly increase the search space. In addition, 
converting constraints into squared penalty terms often entails     
all-to-all connectivity between a problem's variables, which 
can be challenging to realize for hardware architectures with nearest-neighbour connectivity or other geometric locality constraints,
necessitating expensive \textsc{Swap} gate routing strategies 
that lead to deeper circuits.

We demonstrate these shortcomings using integer programs called binary linear problems with inequality constraints. The canonical formulation of such problems is:
\begin{align}\label{eq:linear_problem}
    &\min_{x} q_0^T x \nonumber \\
    &\textrm{subject to} \nonumber\\
    &q_i^T x  \geq c_i \quad \textrm{for} \quad i = 1, ...,m \nonumber \\
    &x \in \{0,1\}^n, 
\end{align}
where $q_0 \in \mathbb{R}^n$, $q_i \in \mathbb{Z}^{n}$, and $c_i \in \mathbb{N}$, for $i = 1,...,m$.
To convert \cref{eq:linear_problem} into a QUBO problem, we incorporate all constraints as quadratic penalties with additional slack variables follows:
\begin{widetext}
\begin{align}\label{eq:qubo_linear}
    &\min \left(q_0^Tx + \sum_{i=1}^m \gamma_i \left( q_i^T x - c_i - W_i \right)^2 \right)\nonumber \\
    &\textrm{subject to} \nonumber\\
    &W_i = \sum_{k=0}^{M_i-2} 2^k y_{k}^{(i)} + R_i y_{M_i-1}^{(i)} \;\textrm{ for} \ i=1, \ldots, m \nonumber\\
    &y^{(i)} \in \{0,1\}^{M_i} \;\textrm{ for}  \ i=1,\ldots, m \nonumber \\
    &x \in \{0,1\}^n.
\end{align}
\end{widetext}
In this problem, we minimize an objective function over binary vectors $x$ and $y^{(i)}$ for $i=1,\ldots,m$. The scalar $\gamma_i > 0$ is a penalty coefficient, and $0 \leq W_i \leq c_i^{\text{max}} :=\max_x \{ q_i^T x - c_i \}$ is an integer slack variable given by a binary expansion in terms of \mbox{$M_i=\lfloor \log_2(c_i^{\text{max}}) \rfloor + 1$} additional auxiliary binary variables $y_{k}^{(i)} \in \{0,1\}$. The constant $R_i$ denotes a remainder which is chosen such that $W_i \leq c_i^{max}$ for all vectors $y^{(i)}$. Note that \cref{eq:qubo_linear} is a quadratic polynomial over binary variables $x_j$ for $j = 1,\ldots n$ and $y_{k}^{(i)}$ for \mbox{$k = 0,\ldots, \lfloor \log_2(c_i^{\text{max}}) \rfloor$} and $i=1, \ldots, m$. If a constraint is satisfied, $q_i^Tx \geq c_i$, then the quadratic penalty is zero, $\gamma_i \left( q_i^T x - c_i - W_i \right)^2=0$. However, if a constraint is violated, $q_i^Tx < c_i$, then we have a non-zero quadratic penalty, \mbox{$\gamma_i \left( q_i^T x - c_i - W_i \right)^2 > 0$}, because $W_i\geq 0$. In the following subsections, we demonstrate various issues resulting from QUBO reformulation through the lens of this problem.

\subsection{\label{sec:qubo_issues_additional_qubits}  Overhead Resulting from  Auxiliary Variables}
From \cref{eq:qubo_linear} it is clear that the reformulation of the original constrained problem \cref{eq:linear_problem} as QUBO has significantly more binary variables. Concretely, the new problem has $n + \sum_{i=1}^m M_i$ variables instead of $n$ variables in the canonical formulation \cref{eq:linear_problem}. This implies a larger search space of dimension $2^{n + \sum_{i=1}^m M_i}$, necessitating  
$\sum_{i=1}^m M_i$ additional qubits for encoding the circuit. Recall that $M_i = \lfloor \log(c_i^{\text{max}}) \rfloor +1$. If $c_i^{\text{max}}$ is large, the overhead associated with the auxiliary qubits may become much larger than $n$, which implies that even relatively small problem instances may become highly challenging when converted to a QUBO problem.

\subsection{\label{sec:qubo_issues_connectivity} Significantly Increased  Connectivity}
Reformulations such as QUBO \cref{eq:qubo_linear}  typically also affect  connectivity requirements. Suppose that for some $i_0 \in \{1,\ldots,m\}$, we have a vector $q_{i_0}$ with no zero components, i.e., $q_{i_0,j} \neq 0$ for all $j=1,\ldots, n$. Then the quadratic penalty
$\gamma_{i_0}(q_{i_0}^Tx - c_i - W_{i_0})^2$ has $\binom{n + M_{i_0}}{2}$ quadratic terms of the form $x_j x_v$, $x_j y_k^{(i_0)}$ and $y_k^{(i_0)} y_u^{(i_0)}$ for $j,v = 1, \ldots, n$ and $k, u = 0,\ldots, \lfloor \log_2(c_{i_0}^{\text{max}}) \rfloor$. That is, every variable is coupled with every other variable, thus requiring all-to-all connectivity between the qubits representing them. This implies that every layer of a QUBO-based quantum circuit will have $\binom{n + M_{i_0}}{2}$ 2-qubit $R_{ZZ}$ gates 
to be applied on each possible pair of qubits. This makes the QAOA circuit computationally demanding due to all-to-all connectivity.

\begin{figure}[h]
    \centering
    \includegraphics[scale=0.7]{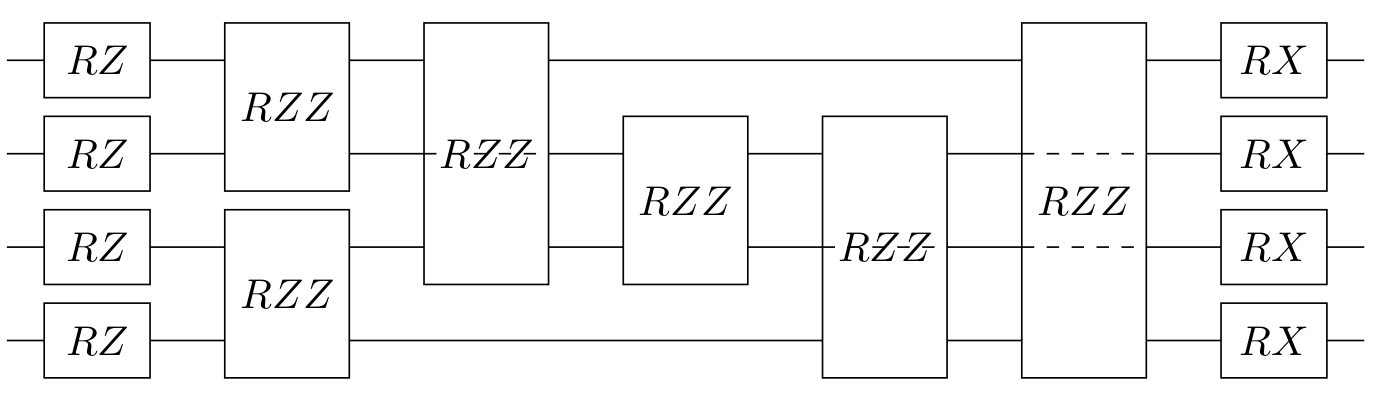}
    \caption{\label{fig:qubo_circuit} Typical quantum circuit structure 
    for a single layer when implementing DAQC for a 4-qubit QUBO problem. Rotation angles are not displayed for clarity. The circuit is all-to-all connected since each qubit interacts with every other qubit through $R_{ZZ}$ gates. The dashed wires plotted over some $R_{ZZ}$ gates denote wires that are not affected by the gate.}
\end{figure}

Moreover, if each circuit layer involves all-to-all connectivity, the circuit depth depends on the number of layers, denoted in this paper by the parameter $p$, and the QUBO problem size \mbox{$N = n + \sum_{i=1}^m M_i$}. In this case, the circuit depth is $O(Np)$. For more details on this point, see~\cref{ap:complexity}. A typical circuit structure for a QUBO problem 
is illustrated in~\cref{fig:qubo_circuit}.

\subsection{Challenges from Optimization Perspective}
From the optimization point of view, each additional variable $y_k^{(i)}$ doubles the search space. Hence, the size of the search space increases exponentially with respect to the number of additional variables. In addition, the energy landscape generally becomes more rugged, because flipping the value of $x_j$ or $y_k^{(i)}$ may lead to high peaks of energy. For example, consider a feasible solution vector ${(x^*,y^*) \in \{0,1\}^{n + \sum_{i=1}^m M_i}}$ such that all penalty terms are zero. If we now flip a single bit $y^{*(i)}_k$ for $k \geq 1$, then the penalty term becomes $\gamma_i \left( q_i^T x^* - c_i - W_i \right)^2 = \gamma_i \left (2^k \right)^2$. This suggests that solutions differing by a single-bit flip result in exponentially high energy peaks, creating a highly rugged optimization landscape. As a result, the associated optimization problem becomes notably more difficult.

\subsection{Concluding Remarks About QUBO}
In this section, we have seen that some CO problems, when converted to QUBO, may require significantly more quantum resources. Problems with inequality constraints require additional auxiliary binary variables, and the number of such variables grows logarithmically with the maximum constraint's value. Quadratic penalties create additional interactions between variables, giving rise to additional 2-qubit $RZZ$ gates. From the optimization perspective, the QUBO formulation has a larger search space, and this makes the problem more difficult for any type of algorithm. Moreover, solutions that are one-bit flip away from the current feasible solution are often infeasible and yield exponentially large penalties.

\section{\label{sec:lagrangian_dual_approach}Lagrangian dual approach}
This study departs from the conventional QUBO approach and investigates the solving of constrained problems using Lagrangian duality theory. We demonstrate the superiority of the proposed approach by comparing it with the QUBO-based approach in the setting of DAQC. 

\subsection{\label{sec:lagrangian_duality}Lagrangian Duality}
We propose a Lagrangian dual DAQC protocol for solving binary combinatorial problems with inequality constraints of the form given in \cref{eq:linear_problem}.
We also note that the proposed approach generalizes easily to binary quadratic problems with quadratic inequality constraints of the following form:
\begin{align}\label{eq:original_quad_problem}
    P := &\min_{x} x^T Q_0 x\nonumber\\
    &\textrm{subject to} \nonumber\\
    &x^T Q_i x  \geq c_i \quad \textrm{for} \quad i = 1, \ldots ,m \nonumber \\
    &x \in \{0,1\}^n, 
\end{align}
where $Q_i \in \mathbb{R}^{n \times n}$ is a symmetric matrix for $i=0,\ldots,m$. The applicability of Lagrangian duality for such problems is investigated in \cite{karimi2017subgradient}. If $Q_i$ is diagonal for $i=0, \ldots, m$, then \cref{eq:original_quad_problem} is equivalent to~\cref{eq:linear_problem} due the relation $x_i^2 = x_i$. Therefore, the formulation in \cref{eq:original_quad_problem} covers a wide range of linear (e.g. \cref{eq:linear_problem}) and quadratic constrained and unconstrained combinatorial problems.

We address the QUBO issues discussed in previous sections by using  Lagrangian duality which is often used in classical optimization. The Lagrangian dual problem corresponding to the primal constrained problem in \cref{eq:original_quad_problem} is
\begin{align}\label{eq:lagrangian}
	D := \max\limits_{\lambda \in \mathbb{R}^m_{+}} \;\min \limits_{x \in \{0,1\}^n} \left(x^TQ_0x + \sum_{i=1}^m \lambda_i \left (c_i - x^TQ_ix\right)\right),
\end{align}
where $\lambda_i \geq 0$ for $i = 1,..., m$ are non-negative Lagrange multipliers. Generally, weak duality holds. That is, if $(x^*, \lambda^*)$ is an optimal dual pair and $D$ is the corresponding objective function value of \cref{eq:lagrangian}, then  $D \leq P$ where $P$ is the optimal objective function value of the primal problem  defined in \cref{eq:original_quad_problem}. For $x^*$ to be an optimal solution to the primal problem [\cref{eq:original_quad_problem}], it must be feasible and satisfy the complementary slackness condition \cite{geoffrion1974lagrangean} below:
\begin{equation}
    \sum_{i=1}^m \lambda_i^* (c_i - x^{*T} Q_i x^*) = 0.
\end{equation}
Whenever this condition holds, the duality gap is zero, that is, $D = P$, and $x^*$ is an optimal solution to the primal problem.

\noindent\textbf{Definition:  $\epsilon$-optimal solution.} If $x^*$ is feasible but the complementary slackness is not satisfied, then \mbox{$\sum_{i=1}^m \lambda_i^* (x^{*T}Q_i x^* - c_i) > 0$}, and we call $x^*$ an \emph{$\epsilon$-optimal solution} to the primal problem \cite{geoffrion1974lagrangean} with $\epsilon$ defined as:
\begin{equation}
    \epsilon:=\sum_{i=1}^m \lambda_i( x^{*T}Q_i x^* - c_i).
\end{equation}
The value  of $\epsilon$ indicates how close $x^*$ is to the optimal solution to the primal problem in terms of its objective function value. That is:
\begin{equation}
    0 \leq x^{*T}Q_0 x^{*} - P \leq \epsilon.
\end{equation}
Furthermore, $\epsilon$ bounds the duality gap:
\begin{equation}
    0 \leq P - D \leq \epsilon.
\end{equation}
Therefore, the Lagrangian dual problem can yield optimal or near-optimal solutions to the primal problem. Unlike QUBO, the formulation \cref{eq:lagrangian} neither requires auxiliary variables nor does it involve squared penalty terms. As a result, it permits significantly more efficient quantum circuits as demonstrated in~\cref{sec:circuit}.

Since our approach is based on adiabatic evolution in time, we may generalize the Lagrange multipliers $\lambda_i$ to time-dependent functions $\lambda_i(t) \geq 0$ for $i=1,...,m$ and $t \in [0, T]$. This allows  enhanced control over the constraint terms during the adiabatic evolution. Hence, we now define the generalized time-dependent Lagrangian relaxation by
\begin{align}\label{eq:lagrangian_time}
	\min \limits_{x \in \{0,1\}^n} \left(x^T Q_0 x + \sum_{i=1}^m \lambda_i(t) \left (c_i - x^T Q_i x \right )\right).
\end{align}
The Lagrangian relaxation is a function of $\{\lambda_i(t)\}_{i=1}^m$. In the context of classical optimization, the generalization of Lagrange multipliers to time-dependent functions is somewhat meaningless. However, in adiabatic quantum computation, which happens through evolution in time, the time-dependent Lagrangian multipliers enable more control over how constraints are introduced over time. A comprehensive discussion of the time-dependent $\lambda_i(t)$ will be undertaken in \cref{sec:lagrangian_mult_schedule}.

A crucial distinction between QUBO reformulations and the Lagrangian dual method is in how their optimal solutions correspond to an optimal solution of the original problem. While the optimal solutions to QUBO reformulations coincide with the optimal solutions of the original primal problem, this is not inherently the case for the Lagrangian dual approach. Due to the principle of weak duality, an optimal solution to the Lagrangian dual problem is usually an $\epsilon$-optimal solution to the primal problem. As a result, classical algorithms for Lagrangian dual problems are not guaranteed to yield the optimal solutions to the primal problem. 

\begin{figure}[t]
\centering
\includegraphics[width=0.4\textwidth]{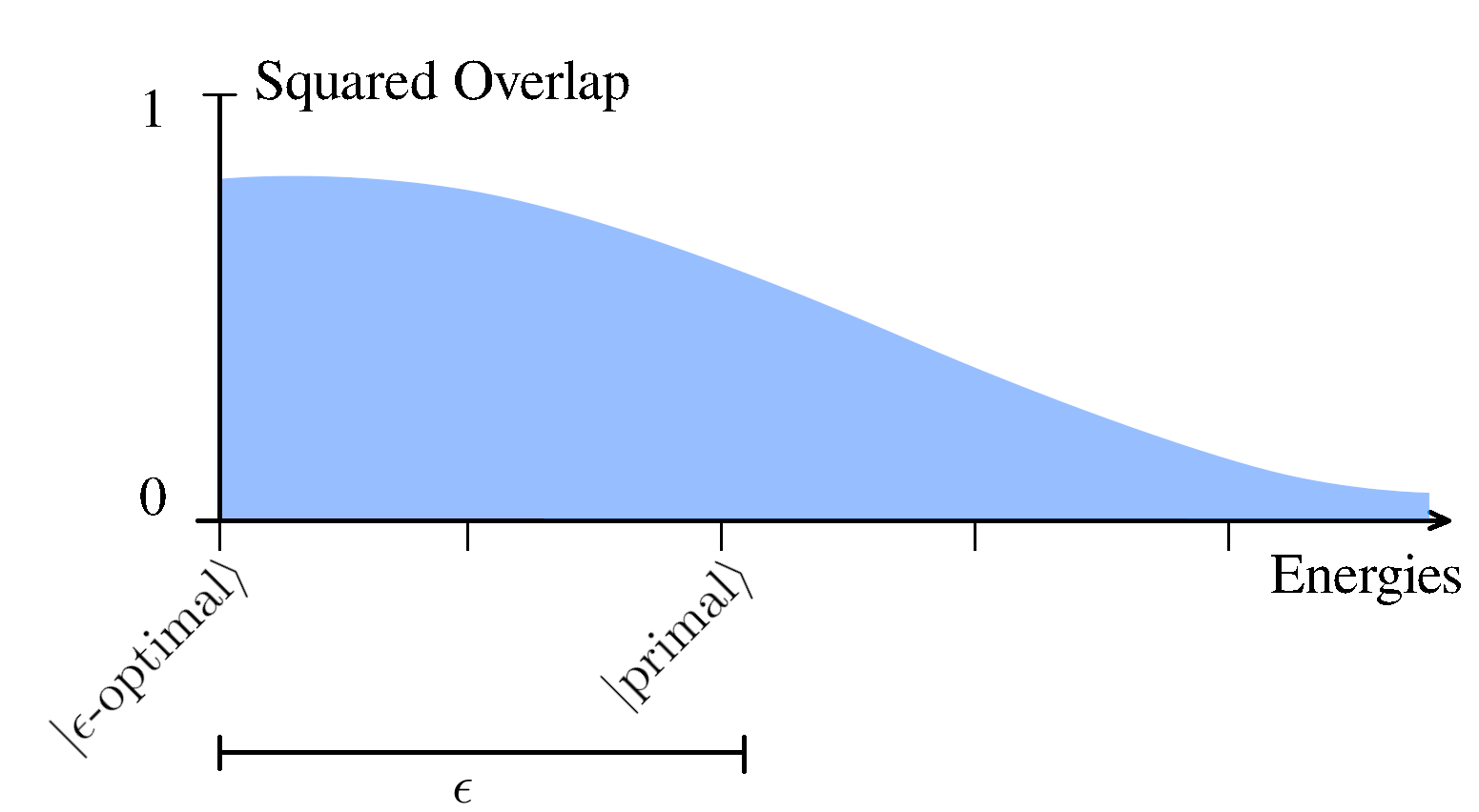}
\caption{
    Schematic illustration of the squared overlap between the state generated by the LD-DAQC circuit and the eigenstates corresponding to the $\epsilon$-optimal and optimal solutions to the primal problem. The $x$-axis represents the energy levels of the Hamiltonian associated with the Lagrangian relaxation, as detailed in \cref{eq:lagrangian_time}. The ground state energy is labeled as $\ket{\epsilon\text{-optimal}}$, while the higher energy level  corresponding to the optimal solution to the primal problem is labeled as $\ket{\text{primal}}$. The $y$-axis measures the amount of squared overlap. Given that the LD-DAQC circuit is engineered to favor the preparation of the lowest energy eigenstate, there is a more significant overlap with eigenstates of lower energies. Because the state of interest $\ket{\text{primal}}$ is at a higher energy level (specified by  the upper bound $\epsilon$), it is easier to reach using shorter evolution times.}
    \label{fig:overlap_illustration}
\end{figure}

Remarkably, in the quantum optimization context, quantum circuits leveraging the Lagrangian duality offer an approach to finding the optimal solutions to the primal problem by optimizing \cref{eq:lagrangian_time}. This is achieved by repeated execution of a pre-tuned \mbox{LD-DAQC} circuit followed by measurements; the associated algorithm is presented in~\cref{sec:circuit}. Our numerical experiments discussed in~\cref{sec:Benchmarking_results} suggest that the optimal solution of the primal problem is often contained in the measured sample. Indeed, the optimal solution to the primal problem is measured frequently enough to outperform the traditional QUBO-based approaches. This observation suggests that quantum states prepared by LD-DAQC circuits have a considerable overlap with the quantum state representing the optimal solution to the primal problem. This can be understood by remembering that the problem Hamiltonian representing \cref{eq:lagrangian_time} is diagonal in the computational basis. Consequently, the ground state of the Hamiltonian corresponds to the $\epsilon$-optimal solution to the primal problem, while the optimal solution to the primal problem [\cref{eq:original_quad_problem}] manifests as an excited eigenstate. It can be shown that $\epsilon$-optimality ensures the energy gap between these two eigenstates is bounded by $\epsilon$. Thus, given a small energy gap and a carefully timed evolution over a finite duration $T$, we anticipate a significant overlap between the state produced by LD-DAQC and the eigenstate representing the optimal solution to the primal problem. If the conditions of strong duality and complementary slackness are met, then with $\epsilon=0$, the $\epsilon$-optimal solution becomes the optimal solution to the primal problem. We illustrate this behaviour in \cref{fig:overlap_illustration} which depicts the squared overlap between the state prepared by LD-DAQC and eigenstates associated with $\epsilon$-optimal and optimal solutions to the primal problem.

We formalize these insights in mathematical generality.
Consider $f(x):\{0,1\}^n \rightarrow \mathbb{R}$ and $g_i(x):\{0,1\}^n \rightarrow \mathbb{R}$ for $i=1,\ldots, m$. The  primal optimization problem is:
\begin{align}\label{eq:general_primal_problem}
    &\min_x f(x) \nonumber \\
    &\text{subject to}\nonumber \\
    &g_i(x) \geq c_i \text{ for } i =1, \ldots , m \nonumber \\
    &x \in \{0,1\}^n.
\end{align}
Then, the Lagrangian dual problem is given by:
\begin{equation}
    \max_{\lambda \in \mathbb{R}^m_{+}} \min_x \left ( f(x) + \sum_{i=1}^m \lambda_i (c_i - g_i(x)) \right).
\end{equation}
Let $L(x)$ denote the Lagrangian function
\begin{equation}\label{eq:general_lagrangian_function}
    L(x) := f(x) + \sum_{i=1}^m \lambda_i^* (c_i - g_i(x)),
\end{equation}
where $\lambda^*>0$ is fixed and optimal.

Let $x_P$ be the optimal solution to the primal problem in \cref{eq:general_primal_problem}, and let $x_D$ be the optimal solution to $L(x)$ in \cref{eq:general_lagrangian_function}. Since $x_D$ is optimal, we have:
\begin{equation}
    L(x_D) \leq L(x_P).
\end{equation}
In the physical context, the value $L(x_D)$ corresponds to the ground state energy of the problem Hamiltonian associated with $L(x)$, and the value $L(x_P)$ corresponds to an elevated energy level, which is easier to reach when the adiabatic transition time is finite. It is straightforward to show  that the gap between $L(x_D)$ and $L(x_P)$ is bounded by $\epsilon$ [\Cref{ap:duality_gap_overlap}]. That is:
\begin{equation}\label{eq:lagragian_gap_epsilon}
    L(x_P) - L(x_D) \leq \epsilon = \sum_{i=1}^m \lambda_i( g_i(x_D) x^* - c_i).
\end{equation}
Therefore, the energies $L(x_P)$ and $L(x_D)$ are $\epsilon$-close, and $\ket{x_D}$ is the $\epsilon$-optimal solution to the primal problem. This proximity affects the overlap between the state $\ket{\psi}$ prepared by LD-DAQC and the eigenstate $\ket{x_P}$. To illustrate this behaviour, we show that, under certain idealized assumptions and simplifying approximations [\Cref{ap:duality_gap_overlap}], the squared overlap depends on $\epsilon$ as follows:
\begin{equation}\label{eq:overlap}
    |\braket{\psi}{x_P}|^2 = \sin^2\left( \frac{A}{Tr} \ln \left (1 + \frac{Tr}{\epsilon} \right)\right),
\end{equation}
where $A = \bra{x_P}H_{\text{init}}\ket{x_D} > 0$, and $r$ is some independent real parameter. Furthermore, it is possible to establish the relation between the parameter $T$ and $\epsilon$. Setting $(A/Tr) \ln(1+Tr/\epsilon) = \pi/2$, we obtain the following relation:
\begin{equation}\label{eq:relation_eps_T}
    \epsilon = \frac{Tr}{e^{\frac{\pi}{2}\frac{Tr}{A}}-1}.
\end{equation}
From the equation above, it is clear that for any $\epsilon>0$, there exist values for $T$, $r$ and $A$ which yield a unit overlap between the state prepared by LD-DAQC and the eigenstate corresponding to the primal problem's solution $x_P$. As mentioned before, in the special case where $\epsilon = 0$, the condition of complementary slackness is satisfied, rendering $\ket{x_D} = \ket{x_P}$ to be the ground state. Therefore, we arrive at the following conclusions. Smaller values of $\epsilon$ require enforcing adiabaticity of the evolution (larger $T$) so that the final state $\ket{\psi}$ has an energy that is $\epsilon$-close to the energy of the ground state $\ket{x_D}$ and hence $\epsilon$-close to the energy of $\ket{x_P}$. On the other hand, for larger values of $\epsilon$, the adiabaticity conditions can be relaxed (smaller values of $T$ to allow transitions to higher energies), as the state of interest $\ket{x_P}$ is energetically further away from the ground state $\ket{x_D}$.

Although \cref{eq:overlap} and \cref{eq:relation_eps_T} offer intuitive and theoretical insights into the LD-DAQC approach, they are based on idealized assumptions. Consequently, these equations are not practical for determining the optimal value of $T$. In real-world scenarios, the identification of the optimal $T$ usually relies on empirical methods, such as hyperparameter tuning. For instance, in dealing with a set of $n$-variable problems, the objective is to find a single value $T$ that yields the best results across the entire problem set.

\subsection{\label{sec:algorithm}Lagrangian Dual and Adiabatic Quantum Computation}
In this section, we use the Lagrangian relaxation in \cref{eq:lagrangian_time} to construct an efficient LD-DAQC circuit. First, we recast \cref{eq:lagrangian_time} as a problem of finding the ground state of a Hamiltonian. This is achieved by substituting each variable $x_j$ with a spin variable $s_j$ using the identity $x_j = (1 - s_j)/2$. Then, each variable $s_j$ is substituted with the Pauli-Z operator $Z_j$ which has eigenvalues $+1$ and $-1$. This substitution relates an abstract binary variable $x_j \in \{0,1\}$ with the eigenvalues of the quantum mechanical spin observable and turns the function in \cref{eq:lagrangian_time} into a Hamiltonian that we denote as $H_{\text{P}}(t)$. Here, the problem Hamiltonian becomes explicitly time-dependent, because $\lambda_i(t)$ depends on time.

We use DAQC to find the ground state of $H_{\text{P}}(t)$. For this purpose, we need an initial Hamiltonian $H_{\text{init}}$ whose ground state is known and can be easily prepared in constant time. A typical initial Hamiltonian is the transverse-field operator $H_{\text{init}}=-\sum_i^n X_i$. The adiabatic evolution is achieved by gradually mixing the initial and problem Hamiltonians according to the relation
\begin{align}\label{eq:transition_hamiltonian}
    H(t) = \left(1 - s(t)\right)H_{\text{init}} + s(t)H_{\text{P}}(t).
\end{align}
In this equation, $s(t) \in [0,1]$ for $t \in [0,T]$ specifies an adiabatic schedule with the requirements $s(0)=0$ and $s(T)=1$. 
The system is initially prepared in the lowest-energy eigenstate (ground state) of $H(0) = H_{\text{init}}$. Then, according to the adiabatic theorem \cite{farhi2000quantum, ambainis2004elementary}, a sufficiently slow transition from $H_{\text{init}}$ to $H_{\text{P}}(t)$ guarantees that the system remains arbitrarily close to
the instantaneous ground state of $H(t)$ and thus evolves into the ground state of $H_{\text{P}}$ at the end of the schedule at $t=T$.

It is important to note that, whenever \cref{eq:lagrangian_time} is linear, such as in the knapsack problem that we study in \cref{sec:numerical}, the objective function does not have quadratic terms of the form $x_ix_j$. In this case, the resulting quantum circuit does not have 2-qubit $RZZ$ gates and, consequently, it cannot create entanglement. Yet, it is known that entanglement is a necessary resource for quantum speedup \cite{jozsa2003role, ding2007review}. To introduce entangling operations, we add quadratic terms without modifying the formulation of the optimization problem. Clearly, any addition of quadratic terms $x_ix_j$ to \cref{eq:lagrangian_time} would change the optimization problem to a different optimization problem; similarly, any addition of operators $Z_iZ_j$ to $H_{\text{P}}(t)$ would change the problem Hamiltonian such that it is no longer equivalent to the problem in \cref{eq:lagrangian_time}. Our solution is to use a different initial Hamiltonian with quadratic terms that do not commute with the $Z_i$ or $Z_i Z_j$ terms in $H_{\text{P}}(t)$. More specifically, we consider the 2-local Hamiltonian referred to as the $ZZXX$ universal Hamiltonian \cite{biamonte2008realizable} which has the form:
\begin{align}\label{eq:general_hamiltonian}
    H_{ZZXX} &= \sum_i h_i X_i + \sum_i \Delta_i Z_i \nonumber \\
    &+ \sum_{i,j} K_{i,j}X_i X_j + \sum_{i,j} J_{i,j} Z_i Z_j.
\end{align}
We identify the terms $\sum_i h_i X_i$ and $\sum_{i,j} K_{i,j}X_i X_j$ with the new initial Hamiltonian, whereas the rest of the terms are used to represent the problem Hamiltonian. Therefore, 2-qubit interactions are introduced through the use of the terms in $\sum_{i,j} K_{i,j} X_i X_j$ of the 2-local $ZZXX$ Hamiltonian. We hypothesize that a specific choice of the coefficients $K_{i,j}$ can introduce correlations between qubits that could potentially yield better performance. Generally, the choice of $K_{i,j}$ can be informed by quantum hardware architecture or by the structure of a CO problem. In this study, we choose the coupling coefficients $K_{i,j}$ such that they form a chain with a periodic boundary condition. Specifically, we use the mixing Hamiltonian 
\begin{equation}\label{eq:mixer}
	H_{\text{init}} = -\sum_{i=1}^n X_i - \sum_{i = 1}^{n} X_i X_{i+1}\,,
\end{equation}
where we define $n+1 := 1$. As shown in~\cref{sec:circuit}, this choice yields a highly parallelizable circuit with circuit depth independent of the problem size $n$.

Note that, if the problem in \cref{eq:original_quad_problem} is linear ($Q_i$ for $i=0,\ldots, m$ are diagonal), the associated Lagrangian dual in \Cref{eq:lagrangian} is also linear. It follows that $H_{\text{P}}$ is a 1-local Hamiltonian and $H_{\text{init}}$ in \cref{eq:mixer} is a 2-local Hamiltonian with a user-provided qubit coupling. It is also possible to set $K_{i,j} = 0$ for all $i,j$. Then the Hamiltonian $H(t)$ becomes a 1-local Hamiltonian, because both $H_{\text{init}}$ and $H_{\text{P}}$ are 1-local Hamiltonians. In this particular case, the simulation of evolution generated by $H(t)$ is in the complexity class~{\bf P}, and the corresponding LD-DAQC circuit can be efficiently simulated classically. Therefore, by appropriately selecting the coefficients $K_{i,j}$, one could devise two distinct heuristics: one that leverages entanglement to achieve a potential quantum speedup, and another one that uses Hamiltonian 1-locality to obtain an efficient classical heuristic.

\section{\label{sec:circuit}LD-DAQC for Linear Problems}
In this section, we construct an efficient quantum circuit which, unlike circuits derived from the QUBO approach, is substantially more parallelizable, does not require auxiliary qubits, and requires only nearest-neighbour qubit connectivity.

To find the optimal solution to \cref{eq:linear_problem}, we start the development of our quantum circuit by formulating the time-dependent Lagrangian relaxation:
\begin{align}\label{eq:new_formulation_general}
   \min \limits_{x \in \{0,1\}^n} \left(q_0^t x + \sum_{i=1}^m \lambda_i(t) \left (c_i - q_i^T x\right) \right).
\end{align}
Generally, $\lambda_i(t)$ can be a parametrized function of time as in \cref{eq:lambda_sch}, and its parameters are determined through hyperparameter tuning as explained in \cref{sec:exp_setup}. In the special case when $\lambda_i(t) \equiv \lambda_i$ is time-independent for all $i$, then this objective function is convex in $\lambda_i$. The optimal $\lambda_i^*$ can be estimated to an arbitrary error $\epsilon > 0$ in $O(\epsilon^{-2})$ iterations using the subgradient method \cite{shor2012minimization}.

Next, we recast \cref{eq:new_formulation_general} as a problem of finding the ground state of a Hamiltonian. By substituting binary variables $x_j$ with spin variables $s_j$ according to the relation $x_j = (1 - s_j)/2$ and introducing \mbox{Pauli-Z} operators $Z_j$ we obtain the problem Hamiltonian $H_{\text{P}}(t)$ which is linear in $Z_j$. Dropping all constant terms in $H_{\text{P}}(t)$ yields
\begin{align}\label{eq:kp_h_problem}
    H_{\text{P}}(t) = \sum_{j=1}^n \left(-q_{0j} + \sum_{i=1}^m \lambda_i(t) q_{ij}\right) Z_j.
\end{align}
Combining $H_{\text{P}}(t)$ with $H_{\text{init}}$ in \cref{eq:mixer} yields the total Hamiltonian
\begin{align}
    H(t) = & \left ( 1 - s(t) \right ) \left (-\sum_{j=1}^n X_i - \sum_{j = 1}^{n} X_j X_{j+1} \right ) \nonumber \\
    + & s(t) \sum_{j=1}^n \left(-q_{0j} + \sum_{i=1}^m \lambda_i(t) q_{ij}\right) Z_j, \ t \in [0, T],
\end{align}
where $s(t)$ is an adiabatic schedule. Empirically, it was found \cite{sankar2021benchmark}  that normalization of $H_{\text{init}}$ and $H_{\text{P}}(t)$ by their respective Frobenius norms yields a well-performing schedule. Therefore, we use the following rescaled and dimensionless total Hamiltonian:
\begin{equation}
    \tilde{H}(t) = (1 - s(t))\frac{H_{\text{init}}}{||H_{\text{init}}||} + s(t)\frac{H_{\text{P}}(t)}{||H_{\text{P}}(t)||},
\end{equation}
with $H_{\text{init}}$ and $H_{\text{P}}(t)$ given by  \cref{eq:mixer} and \cref{eq:kp_h_problem}, respectively. Given the initial ground state $\ket{\psi(0)}=\ket{+}^{\otimes n}$ of $H_{\text{init}}$, the evolved state $\ket{\psi(t)}$ can be expressed as $\ket{\psi(t)} = U(t,0)  \ket{\psi(0)}$, where $U(t,0)$ is the time evolution operator which satisfies the dimensionless Schr\"{o}dinger equation,
$ \mathrm{i} dU(t,0)/dt = \tilde{H}(t)U(t,0)$, with a rescaled and dimensionless time variable.

In order to construct the LD-DAQC circuit,  we subdivide the time interval $[0, T]$ into $p$ subintervals of length ${\Delta t = T / p}$. Then the evolution operator is given by
\begin{align}\label{eq:exact_solution}
    U(T,0) = \lim_{p \rightarrow \infty} \prod^p_{k=0} \exp \left \{-\mathrm{i} \Delta t \tilde{H}(\Delta t k) \right \}.
\end{align}
To approximate \cref{eq:exact_solution}, we choose $p$ to be finite and apply the established first-order Trotterization formula. This yields a $p$-layered LD-DAQC circuit, 
\begin{align}\label{eq:lag_approx_U}
    \hat{U}(T, 0) = \prod^{p}_{k=1} \exp \left \{ -\mathrm{i} \gamma_k H_{\text{init}} \right \} \exp \left \{ -\mathrm{i} \beta_k H_{\text{P}}(\Delta t k) \right \},
\end{align}
with $\gamma_k, \beta_k$ defined as
\begin{align}\label{eq:coeffs}
    \gamma_k &= \frac{\left ( 1 - s(k \Delta t) \right ) \Delta t}{||H_{\text{init}}||}, \nonumber \\
    \beta_k &= \frac{s(k \Delta t ) \Delta t}{||H_{\text{P}}(k \Delta t)||}.
\end{align}

\begin{figure}[t]
    \centering
    \includegraphics[scale=0.85]{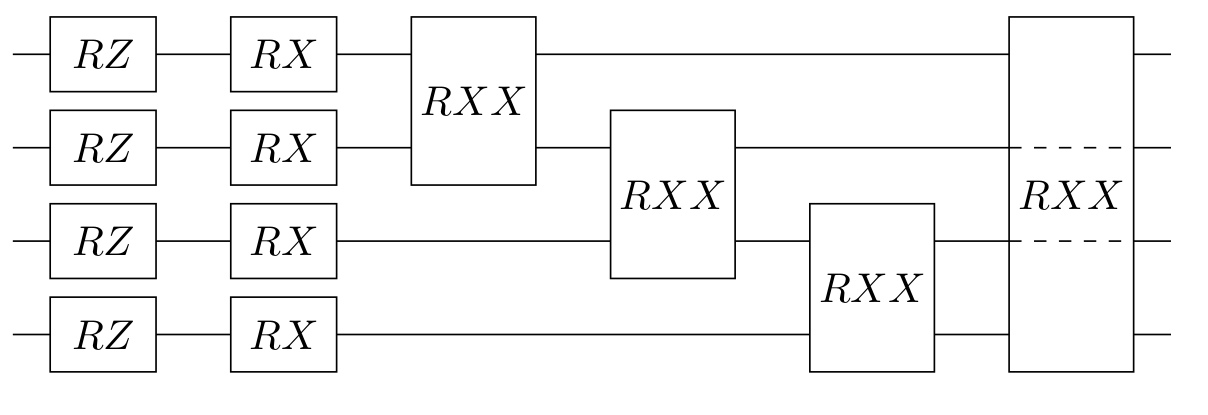}
    \caption{\label{fig:daqc_circuit} The $k$-th layer of the 4-qubit LD-DAQC. Rotation angles are not displayed for clarity.}
\end{figure}

We now examine the structure of the resulting circuit $\hat{U}(T, 0)$ given in \cref{eq:lag_approx_U}. We will show that the circuit depth complexity depends only on the parameter $p$ and it does not explicitly depend on the problem size $n$. To see this, for any layer $k \in \{ 1,...,p\}$, consider the unitary matrix given by the initial Hamiltonian. Due to the commutativity of $X_i$ and $X_iX_j$ the matrix can be written as
\begin{align}\label{eq:gates_init}
    \exp \left \{ -\mathrm{i} \gamma_k H_{\text{\text{init}}} \right \} =
    \prod_{j=1}^n \exp \left \{ \mathrm{i} \gamma_k X_j X_{j+1} \right \} \prod_{j=1}^n \exp \left \{ \mathrm{i} \gamma_k X_j \right \}.
\end{align}
The unitary matrix given by the problem Hamiltonian can be written as
\begin{align}\label{eq:gates_problem}
    \exp & \left \{ -\mathrm{i} \beta_k H_{\text{P}}(\Delta t k) \right \} = \nonumber \\
    &\prod_{j=1}^n \exp \left \{ -\mathrm{i} \beta_k \left(-q_{0j} + \sum_{i=1}^m \lambda_i(\Delta t k) q_{ij}\right) Z_j \right \}.
\end{align}
We note that the Lagrange multiplier $\lambda_i(\Delta t k)$ in \cref{eq:gates_problem} contributes to the angle of rotation and depends on the time step $\Delta t k$. For any $k$, the unitary matrices in \cref{eq:gates_problem} can be represented by 1-qubit $RZ$ gates. Conversly, the unitary matrices given in \cref{eq:gates_init} can be represented by 1-qubit $RX$ gates and 2-qubit $RXX$ gates. Hence, the LD-DAQC circuit $\hat{U}(T, 0)$ can be expressed as follows:
\begin{align}\label{eq:lagrangian_general_circuit}
    \hat{U}(T, 0) = &\prod_{k=1}^p
    \prod_{j=1}^n RXX_{j,j+1}(-2 \gamma_k) \nonumber \\
    \times&\prod_{j=1}^n RX_j(-2 \gamma_k) \nonumber \\
    \times&\prod_{j=1}^n RZ_j\left( 2 \beta_k \left( -q_{0j} + \sum_{i=1}^m \lambda_i(\Delta t k) q_{ij} \right)\right).
\end{align}
The structure of the $k$-th layer of $ \hat{U}(T, 0)$ is illustrated in Fig.~\ref{fig:daqc_circuit}.
The resulting circuit is extremely parallelizable and hence has a short circuit depth. All $RZ$ gates can be applied in a single time step. Similarly, all $RX$ gates can be applied in a single time step. All $RXX$ gates can be applied in 2 time steps if $n$ is even,  and in 3 time steps if $n$ is odd. Hence, the circuit depth complexity is $O(p)$; that is, the circuit depth depends only on the parameter $p$, but it does not depend on the problem size;  
for further discussion, see~\cref{ap:complexity}. Moreover, due to the choice of the coefficients $K_{i,j}$ in $H_{\text{\text{init}}}$ in \cref{eq:mixer}, the resulting circuit features only nearest-neighbour qubit connectivity. Finally, this circuit construction does not involve auxiliary qubits.

\subsection{\label{sec:daqc_schedule}DAQC Scheduling}
In this section, we give an explicit description of the adiabatic schedule $s(t)$ used in this study. Ideally, we would  like to have a single adiabatic schedule that maximizes the probability of observing optimal solutions to multiple problem instances of the same problem class. A problem class can be, for example, an $n$-variable knapsack problem or an $n$-variable travelling salesperson problem.

We adapt the QA scheduling presented in \cite{roland2002quantum}. For better control of the curvature of the schedule function, we use a cubic polynomial approximation proposed in \cite{sankar2021benchmark}. The general form is given as follows
\begin{align}\label{eq:schedule}
    s(t; a, T) = \frac{t}{T} + a \cdot \frac{t}{T} \left (\frac{t}{T}-\frac{1}{2} \right ) \left (\frac{t}{T}-1 \right ).
\end{align}
The parameters $a$ and $T$ define the slope and the evolution time, respectively. The goal is to determine optimal parameters $a^*$ and $T^*$ that maximize the probability of observing an optimal solution for similar problem instance of the same size. This is usually achieved by using a ``training'' set of problem instances of size $n$ and performing a Random Search optimization \cite{bergstra2012random} on $a$ and $T$ to directly minimize the median of \mbox{\em time-to-solution} (TTS) metrics (to be defined in \cref{sec:metrics}) for all problems in the ``training'' set. Once $a^*$ and $T^*$ are determined, we can reuse $s(t; a^*, T^*)$ across all similar problem instances of size $n$ and achieve a high probability of measuring the optimal solution.

\subsection{Lagrangian Multiplier Scheduling \label{sec:lagrangian_mult_schedule}}
In this section we discuss the generalized Lagrange multipliers $\lambda_i(t)$ for $i=1,...,m$ and $t \in [0, T]$. First, we note that constant $\lambda_i$ is a special case of $\lambda_i(t)$. Most combinatorial problems have no notion of time. That is, the problem formulation and its solution are time-independent. Therefore, it makes sense to use Lagrange multipliers that are also constant. However, when using the DAQC approach, the solution to a problem is gradually obtained during the adiabatic evolution. The inherent time-dependence of the adiabatic process can be used to control the strength of the Lagrange multipliers. For example, it is possible to delay the introduction of some difficult constraints allowing the process to start with an easier problem. Alternatively, soft constraints, which are preferred but do not necessarily have to be fulfilled, can also be scheduled to appear towards the end of the evolution. In this study, we define the time-dependent Lagrangian multipliers in terms of the modified schedule in \cref{eq:schedule}:
\begin{align}\label{eq:lambda_sch}
    \lambda_i(t; o_i, \gamma'_i, a_i) := \gamma'_i \cdot s(t - o_i; a_i, T) \cdot \mathbb{1}_{o_i < t}(t) \nonumber \\
    \textrm{ for } \ i = 1, \ldots, m. 
\end{align}
In the equation above, $\gamma'_i > 0$ is the weight of the schedule,  $o_i \in [-T, T]$ is the time offset, $a_i$ is a schedule slope coefficient, and $\mathbb{1}_{o_i < t}(t)$ is an indicator function such that $\mathbb{1}_{o_i < t}(t) = 1$ for $o_i < t$ and $\mathbb{1}_{o_i < t}(t) = 0$ otherwise. Setting $o_i = T/2$ introduces the constraint $i$ in the middle of the adiabatic process, whereas $o_i = T$ is equivalent to ignoring the constraint completely as $\lambda_i(t; o_i, \gamma'_i, a_i) = 0$ for all $t \in [0, T]$.

\section{\label{sec:numerical}Numerical experiments}
This section presents a well-known linear problem that is a special case of the general formulation in \cref{eq:linear_problem}, and compares numerical results obtained by QUBO-DAQC and LD-DAQC circuits for that problem. We will see that the proposed LD-DAQC significantly outperforms the canonical QUBO-based circuit.

\subsection{Knapsack Problem}\label{sec:KP problem}
For our benchmarking study, we use a well-known optimization problem called 1D 0--1 {\em knapsack problem} (KP) \cite{kellerer2004introduction} as the reference problem for comparison. A KP is an NP-complete linear problem that can be formulated as a special case of \cref{eq:linear_problem} in which $m=1$. In combinatorial optimization, a KP requires selecting the most valuable items to fit into a knapsack without exceeding its capacity. This problem has diverse industrial applications, including resource allocation, scheduling, and inventory management.

Our choice of this problem is based on several considerations. First, a KP is a typical representative of NP-complete constrained binary linear problems. Second, as we will show below, trying to solve the KP in the quantum setting by the means of adiabatic or variational methods with the canonical QUBO reformulation [\cref{eq:qubo_linear}] results in a prohibitively costly circuit which requires all-to-all connectivity and additional auxiliary qubits whose number grows logarithmically with the constraint bound. Therefore, a KP is a good representative of CO problems that challenge currently available quantum optimization approaches.  Developing a quantum protocol that successfully solves the KP means that the majority of other binary linear programs with inequality and equality constraints can also be successfully solved by the same protocol.  

We derive a KP from the general canonical formulation given in \cref{eq:linear_problem}. Suppose we have a KP with $n$ items. Let $v_j, w_j \in \mathbb{N}$ denote the $j$th item's value and weight respectively and $c \in \mathbb{N}$ is the knapsack capacity. By setting the entries $q_{0,j} = -v_j$, $q_{1,j} = -w_j$ for all $j=1,\dots,n$ and $c_1 = -c$  we obtain the canonical formulation of the $n$-variable KP
\begin{align}
    & \textrm{max} \sum_{j=1}^n v_j x_j \nonumber\\
    & \textrm{subject to} \nonumber\\ 
    & \sum_{j=1}^n w_j x_j \leq c \nonumber \\
    & x \in \{0,1\}^n. \label{eq:kp-problem}
\end{align}

\subsection{LD-DAQC Circuit}\label{sec:lagrangian_relaxation}
We start the development of our LD-DAQC by approximating the problem in \cref{eq:kp-problem} using the Lagrangian dual given below:
\begin{align}\label{eq:new_formulation_kp}
	\min_{x \in \{0,1\}^n} \left(-\sum_{j=1}^n v_j x_j + \lambda(t) \left( \sum_{j=1}^n w_j x_j - c \right)\right).
\end{align}
As before, if we let $\lambda(t) \equiv \lambda$ to be independent of time, then the problem \cref{eq:new_formulation_general} is convex in $\lambda$. Note that \cref{eq:new_formulation_kp} is a particular case of the general formulation in \cref{eq:new_formulation_general}. From the general LD-DAQC given in \cref{eq:lagrangian_general_circuit} it is straightforward to deduce the structure of the LD-DAQC for the KP; we set $m=1$, $q_{0j} = -v_j$ and $q_{1j} = w_j$. This gives the LD-DAQC for the KP.

\subsection{QUBO-DAQC Circuit}
To derive a QUBO circuit for KP, we first reformulate the KP in \cref{eq:kp-problem} as a QUBO problem. For this we use the general QUBO formulation given in \cref{eq:qubo_linear}. In this case, we get:
\begin{align}
	&\min_{x, y} \left(-\sum_{j=1}^n v_j x_j + \gamma \left( \sum_{j=1}^n w_j x_j - W \right )^2 \right)\nonumber\\
	&\textrm{subject to} \nonumber\\
	&W = \sum_{k=0}^{\lfloor \log_2(c) \rfloor - 1} 2^k y_k + (c + 1 - 2^{\lfloor \log_2(c) \rfloor})y_{\lfloor \log_2(c) \rfloor}\nonumber \\
	&y \in \{0,1\}^{\lfloor \log_2(c) \rfloor + 1} \nonumber \\
	&x \in \{0,1\}^n, 
 \label{eq:qubo}
\end{align}
where $\gamma$ is a penalty multiplier and $c$ is the constraint inequality bound. We note that the number of additional auxiliary binary variables $y_k$ grows as $O(\log_2c)$. Since $w_j > 0$ for $j=1,\ldots,n$, from the discussion in Section~\ref{sec:qubo_issues_connectivity}, it is clear that the squared penalty term introduces $\binom{n + \log_2(c)+1}{2}$ pairwise interactions between every qubit in each layer of the circuit. That is, each layer of the circuit will have $\binom{n + \log_2(c)+1}{2}$ 2-qubit $RZZ$ gates. Therefore, each layer of the circuit is all-to-all connected. Due to this, the circuit is poorly parallelizable and has circuit depth $O\left (p(n + \log_2c) \right)$ -- for further information see Appendix~\ref{ap:complexity}. Therefore, the circuit depth depends on the problem size $n$, the value of the constraint bound $c$ and the number of layers $p$. For example, if the constraint bound $c \geq 2^{2n}$, then the number of auxiliary qubits is at least double of the problem size $n$. This means even relatively small knapsack problems may be very challenging for quantum optimization.

One might attempt to decrease the number of auxiliary qubits by scaling down the coefficients $w_j$ and the constraint bound $c$ by dividing both sides of the inequality constraint in \cref{eq:kp-problem} by some constant $z \in \mathbb{N}$. However, this often leads to fractional values, i.e., $\sum_{j=1}^n w_j x_j / z \in \mathbb{Q} \setminus \mathbb{Z}$ for some $x \in \{0,1\}^n$. This implies that the slack variable $W$ must be able to approximate the rational $\sum_{j=1}^n w_j x_j / z$ for all feasible $x$. Hence, the binary expansion of $W$ must include additional auxiliary binary variables that represent a fractional part. Therefore, in the presence of inequality constraints, it is practically impossible not to have auxiliary qubits.

\section{\label{sec:exp_setup}Experimental setup and Results}
In this section, we discuss the computational setup used in the experiment, present procedures used to generate multiple KP datasets, and various metrics used for our analysis. We compare the performance of QUBO-DAQC and LD-DAQC and demonstrate the superiority of the latter.

\subsection{Hardware Setup}
The experiment involved very large LD-DAQC and QUBO-DAQC circuits, some of which had 100 layers, 20 qubits and all-to-all connectivity. Executing such deep circuits on a quantum computer would require fault-tolerant architectures whose development is still in its infancy. Therefore, all quantum computations were simulated on classical hardware. This was a large-scale experiment which pushed quantum simulation to its limits. To compute the experimental results presented in this paper, we used over 240 hours of 64-core, 32-core and 16-core CPUs and over 200 GB of RAM on the Google Cloud Platform (GCP). Especially demanding QUBO-based quantum circuits whose single execution required at least two weeks of CPU time on the GCP were executed on Nvidia's A100 High-Performance Computing server, reducing the computation time from weeks to several hours. In this experiment, more than 2,000 optimization problem instances were generated and solved using quantum simulation. The correctness of simulation results was verified with Google’s classical optimization solvers: \texttt{OR-Tools Linear Optimization Solver} and \texttt{OR-Tools Brute Force Solver} \cite{cpsatlp}.

\subsection{\label{sec:data_gen}Generation of Random Instances}
We create two different types of \emph{test} supersets of KPs:
\begin{enumerate}
    \item The purpose of the Superset 1 is to examine performance scaling with respect to the problem size $n$. Thus, for each $n=5,6, \ldots,15$ we generate 100 KP instances with integer coefficients $v_i, w_i$ distributed according to the uniform probability distribution with bounds $1$ and $10$. That is, integers $v_i, w_i \sim U(1,10)$. The capacity $c$ is a function of random variables $w_i$ and it is defined as $c = \lfloor \sum_i w_i/2 \rfloor$. 
    \item The purpose of the Superset 2 is to examine performance scaling with respect to the scale of the coefficients $v_i$ and $w_i$ while $n$ is fixed. For a fixed ${n=11}$ and each $C = 10, 20, 30, \ldots, 100$ we generate 100 problem instances such that ${v_i, w_i \sim U(1, C)}$. Hence, there are 10 datasets each containing instances of size $n=11$ but with coefficients $v_i, w_i$ whose range is incremented by 10 units in each dataset.
\end{enumerate}

\subsection{\label{sec:hyperparam_training}Hyperparameter tuning}
We developed two additional \emph{training} supersets by employing the identical setup as previously outlined. These supersets were dedicated solely to hyperparameter tuning. We recall that the DAQC circuits are governed by multiple dimensionless parameters which significantly influence both the performance and circuit complexity:
\begin{itemize}
    \item The hyperparameter $T$ controls the evolution duration. It is important to identify the smallest value of $T$ that satisfies the adiabaticity requirements. Typically, larger values of $T$ necessitate more accurate approximations (more circuit layers or Trotter steps). Less accurate approximations lead to greater errors and thus reduce the success probability of reaching the desired state. Therefore, the hyperparameter $T$ balances the adiabaticity of the evolution and computational efficiency.
    \item The hyperparameter $p$ specifies the number of layers in the DAQC circuit, where higher values of $p$ result in better approximations at the cost of increasing the circuit depth.
    \item The parameter $a$ which determines the curvature of the schedule $s(t)$ in \cref{eq:schedule} and parameters $\gamma'_1$, $o_1$ and $a_1$ specify the Lagrange multiplier scheduling function in \cref{eq:lambda_sch}. These parameters do not affect the complexity of the circuit.
\end{itemize}
To identify the optimal parameters for each $n$, we employ Random Search~\cite{bergstra2012random} across the hyperparameter space, aiming to maximize the probability of sampling an optimal solution of a KP.


\subsection{\label{sec:metrics}Circuit Scaling: LD-DAQC vs.\ QUBO-DAQC}
This section examines the scaling of LD-DAQC and QUBO-DAQC circuits in terms of the \mbox{\em time-to-solution} (TTS) metric. The analysis reveals that LD-DAQC has better scaling and efficiency in TTS than QUBO-DAQC, particularly as problem size and constraints grow.

In order to evaluate the performance of QUBO-DAQC and LD-DAQC on the test datasets we use the $R_{99}$ and the TTS metrics. Let a parameter sequence $\theta = \{(\gamma_k, \beta_k)\}_{k=1}^{p}$ defined in \cref{eq:coeffs} be given. Then $R_{99}(\theta)$ is the number of \lq\lq{}shots\rq\rq{} (executions of the quantum circuit followed by measurements) that must be performed
to ensure a $99\%$ probability of successfully observing the desired state of the problem Hamiltonian $H_{\text{P}}$ for a given circuit parametrization $\theta$ at least once. It is defined as (see, e.g., \cite{sankar2021benchmark})
\begin{align}\label{eq:r99}
    R_{99}(\theta) := \frac{\log(0.01)}{\log(1 - P(\theta))},
\end{align}
where $P(\theta)$ denotes the probability of measuring the optimal solution given the parametrization $\theta$. The TTS denotes the expected computation time required to find an optimal solution for a particular problem instance with $99\%$ confidence. It is defined by
\begin{align}
    \textrm{TTS}(\theta) = R_{99}(\theta)\cdot t_{\text{ss}}\,,
\end{align}
where $t_{\text{ss}}$ denotes the runtime for a single shot. 

In order to estimate $t_{\text{ss}}$, we make the following assumptions:
\begin{enumerate}
    \item A quantum processor performs any single-qubit and two-qubit gate operations in 10 and 20 nanoseconds, respectively; this is in accordance with current state-of-the-art realizations of superconducting quantum hardware technologies, see e.g.\ \cite{suchara2013comparing}.
    \item Gate operations may be performed simultaneously if they do not act on the same qubit.
    \item All components of the circuit are noise-free; hence there is no overhead for quantum error correction or fault-tolerant quantum computation.
\end{enumerate}
Then, as shown in~\cref{ap:complexity}, a $p$-layered LD-DAQC circuit has the following runtime dependence on $p$ and $n$:
\begin{align}
    t_{\text{ss}} = \begin{cases}
    	50p\,, &\;\text{if }\, n \ \textrm{ is even} \\
    	70p\,, & \;\text{if }\,n \ \textrm{ is odd}\;.
    \end{cases}
\end{align}
Note that, for LD-DAQC circuits, $t_{\text{ss}}$ does not explicitly scale with the problem size $n$. In contrast, a $p$-layered \mbox{QUBO-DAQC} circuit has the following runtime dependence on $p$, the problem size $n$,  and the constraint bound $c$:
\begin{align}\label{eq:qubo_runtime}
    t_{\text{ss}} = \begin{cases}
    	20p \left (n + \log_2 c \right), &\text{if }\, n + \log_2 c  \textrm{ is even} \\
    	20p \left (n + \log_2 c + 1 \right), &\text{if }\, n + \log_2c  \textrm{ is odd}\;.
    \end{cases}
\end{align}

Consequently, the TTS scalings for the LD-DAQC and QUBO-DAQC circuits are as follows:
\begin{align}
    \textrm{TTS}^{\text{L}} &= O\left(R_{99}^{\text{L}} \cdot p \right),\\
    \textrm{TTS}^{\text{Q}} &= O\left(R_{99}^{\text{Q}} \cdot p \cdot (n + \log_2 c)\right),
\end{align}
where $\textrm{TTS}^{\text{L}}$ and $\textrm{TTS}^{\text{Q}}$ represent the time to solution for the LD-DAQC and QUBO-DAQC circuits, respectively. It becomes evident that, unless $R_{99}^{\text{L}}$ for the LD-DAQC circuit substantially exceeds $R_{99}^{\text{Q}}$ for the QUBO-DAQC circuit, the LD-DAQC circuit offers a significantly more efficient time to solution. The following section shows that the LD-DAQC circuits require quadratically fewer layers than the QUBO-DAQC circuits to achieve a comparable value for $R_{99}$, i.e., $R_{99}^{\text{Q}} \approx R_{99}^{\text{L}}$.

\subsection{\label{sec:Benchmarking_results} Experiment Results}
In this section, we analyze the optimal resource requirements in terms of the number of layers $p$, evolution time $T$, and the number of qubits for both circuit types to achieve comparable $R_{99}$ values across all problem sizes.
The numerical experiments demonstrate that the LD-DAQC circuits require quadratically fewer resources than QUBO-DAQC circuits while attaining comparable $R_{99}$ values for both approaches. Importantly, we show that $R_{99}$ for LD-DAQC circuits remains unaffected by the expansion in the range of problem coefficients $w_i$ and $v_i$. In contrast, QUBO-DAQC circuits display a polylogarithmic growth in $R_{99}$ when the range of coefficients $w_i$ is increased.

Therefore, the advantage of the LD-DAQC circuits over QUBO-DAQC circuits can be quantified by the ratio of $\textrm{TTS}^{\text{L}}$ to $\textrm{TTS}^{\text{Q}}$. Assume $p$ represents the smallest number of layers required for the QUBO-DAQC circuit so that $R_{99}^{\text{L}} \approx R_{99}^{\text{Q}}$. In this case, the corresponding number of layers needed for the LD-DAQC circuit is $\sqrt{p}$. Hence, we establish the following relationship:
\begin{equation}
    \frac{\textrm{TTS}^{\text{L}}}{\textrm{TTS}^{\text{Q}}} = O\left(\frac{1}{\sqrt p \cdot (n + \log_2 c)}\right).
\end{equation}

\subsubsection{Linear Resouce Scaling of LD-DAQC }
\begin{figure}[h!]
    \centering
    \includegraphics[scale=0.26]{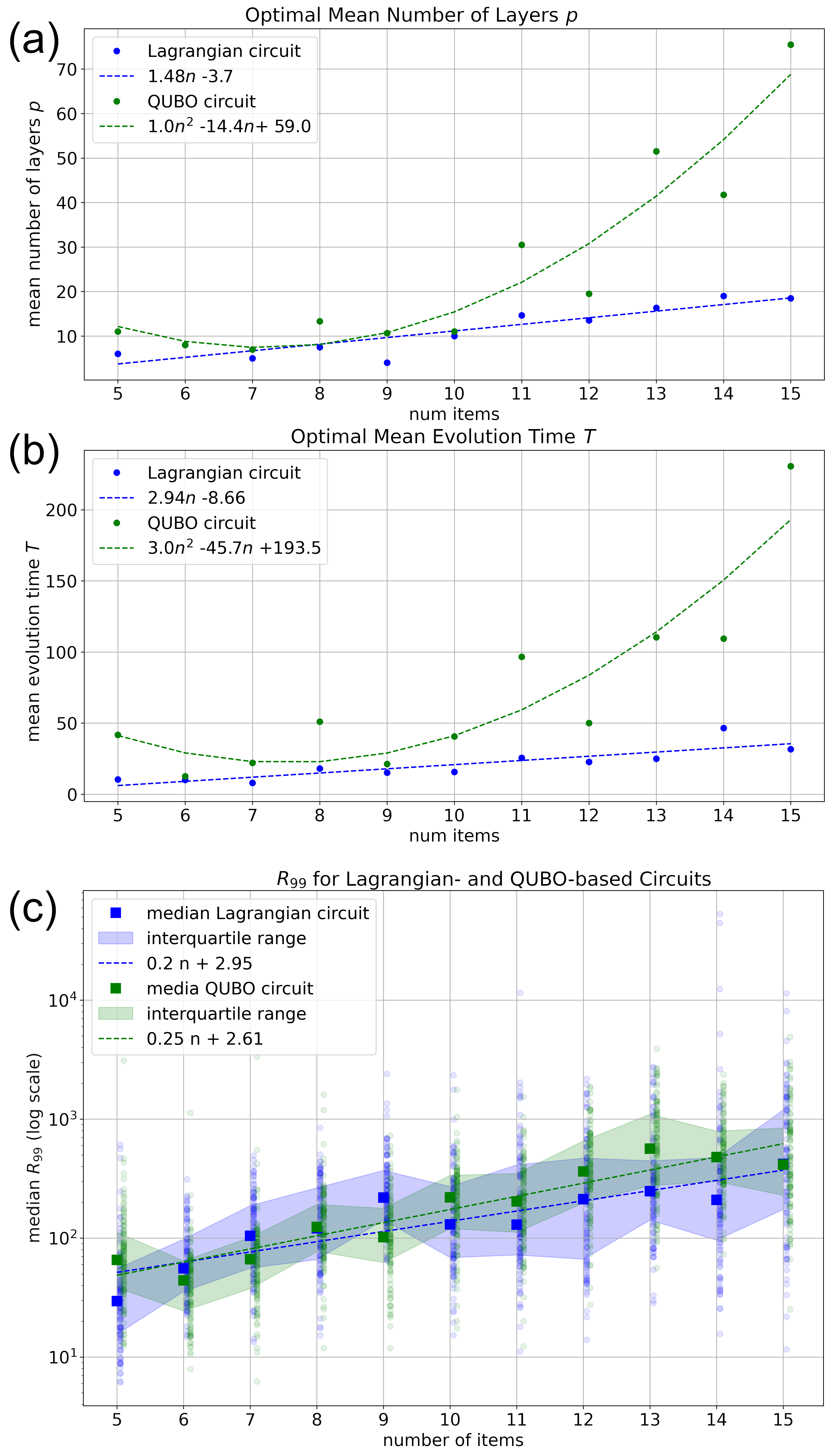}
    \caption{\label{fig:optimal_params_r99} Scaling of the optimal number of layers $p$, optimal evolution time $T$, and $R_{99}$, with respect to the problem size (number of items) for LD-DAQC (blue) and QUBO-DAQC (green) circuits. (a) The optimal mean depth $p$ scales linearly for LD-DAQC, whereas it exhibits a quadratic growth for the QUBO-DAQC circuit. (b) The optimal mean evolution time $T$ also scales linearly for LD-DAQC, whereas it exhibits a quadratic growth for the QUBO-DAQC circuit. 
    (c) $R_{99}$ for LD-DAQC and QUBO-DAQC circuits assuming optimal parameters $p$ and $T$, respectively. Given optimal parameters for each approach, the corresponding circuits have a comparable $R_{99}$ exhibiting  an exponential scaling with respect to the problem size.}
\end{figure}

\cref{fig:optimal_params_r99} displays the scaling of the optimal mean number of layers $p$ and the optimal evolution time $T$ with respect to the problem size. The LD-DAQC circuits have a linear scaling, whereas the QUBO-DAQC circuits have a quadratic scaling for both $p$ and $T$.
 This difference becomes even more critical if we take into account the overhead costs associated with quantum error correction (QEC) and the realization of fault-tolerant quantum computation (FTQC) schemes that become necessary for deep circuits implemented on noisy real-world devices. Deeper circuits are much more vulnerable to noise, necessitating greater QEC and FTQC overhead costs, which have been neglected in our analysis because we aimed at comparing optimistic lower bounds for the two approaches.
\subsubsection{Invariance of LD-DAQC to Coefficient Range Increase}
\begin{figure}
    \centering
    \includegraphics[scale=0.27]{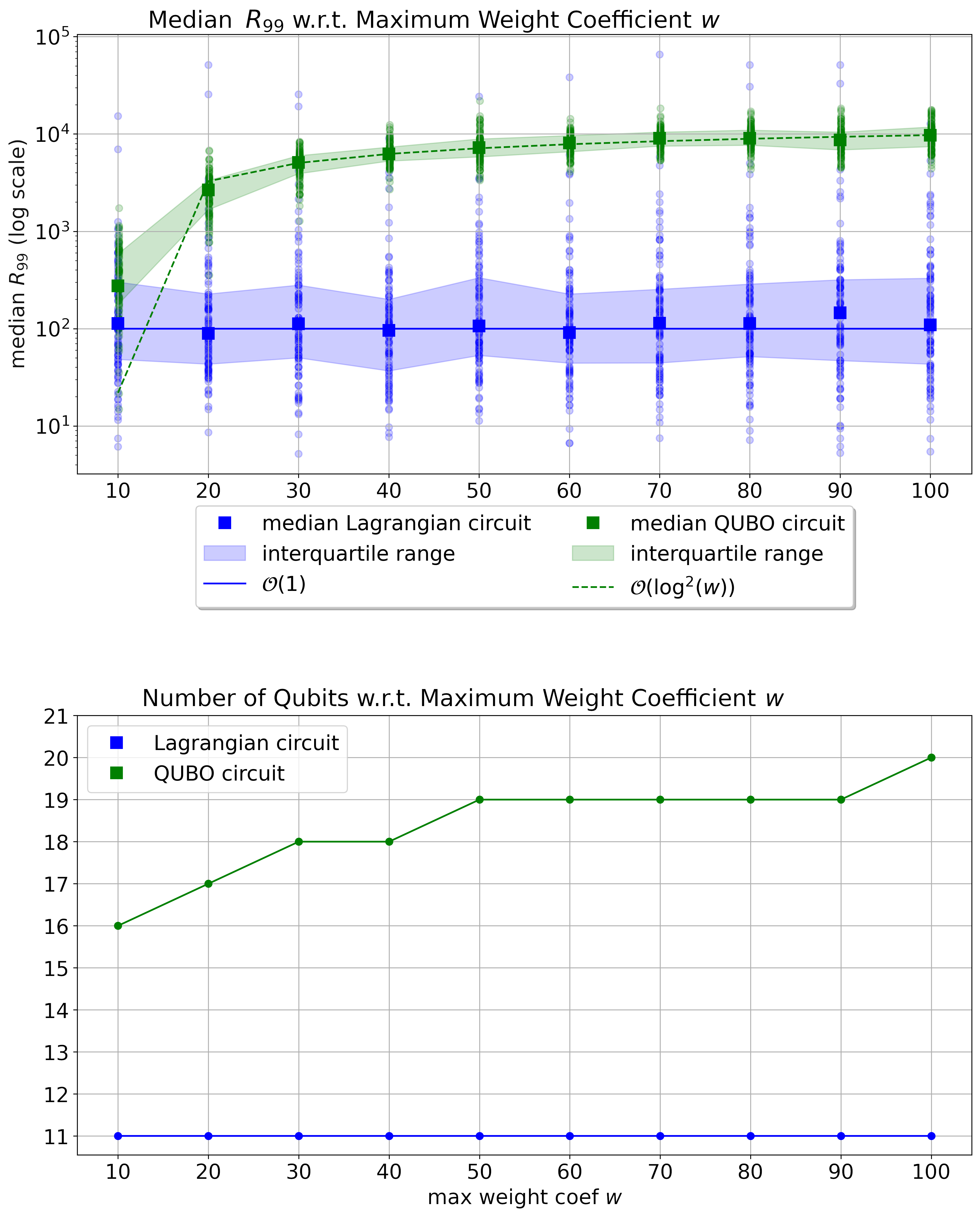}
    \caption{\label{fig:r99_wrt_weight} Scaling of $R_{99}$ and the number of qubits 
    with respect to the maximum weight coefficient (knapsack capacity bound) for LD-DAQC circuits (blue) and QUBO-DAQC circuits (green). A single fixed parametrization was used for all datasets, that is,  $p$, $T$ and $s(t)$ are fixed across all datasets.}
\end{figure}

Recall that whenever a combinatorial optimization problem includes inequality constraints, additional qubits are required to account for encoding the slack variables. It has been established that the count of these additional qubits scales logarithmically with the constraint bounds. Employing Superset 2, we examine the performance scaling on fixed-size KP datasets ($n=11$), but with an increasing value range for the coefficients $v_i$ and $w_i$. In \cref{fig:r99_wrt_weight}, we illustrate the scaling behaviour of $R_{99}$ for both circuit categories, each maintaining fixed parametrization and layer counts across all datasets. It is clear that for LD-DAQC circuits, $R_{99}$ remains invariant to the increase of the range of $v_i$ and $w_i$, exhibiting consistent $R_{99}$ throughout all datasets in Superset 2. In contrast, QUBO-DAQC circuits manifest a quadratically logarithmic growth in $R_{99}$. This observed scaling can be attributed to the logarithmic increase in auxiliary variables intrinsic to QUBO reformulations of problems containing inequality constraints.

The observed scaling behavior in relation to the expanding range of coefficients $v_i$ and $w_i$ underscores a pivotal observation about the TTS scaling of LD-DAQC and QUBO-DAQC circuits. Given the quadratically logarithmic growth of QUBO-DAQC in $R_{99}$ with respect to the increasing range of $w_i$, as shown in \cref{fig:r99_wrt_weight}, we conclude that the TTS disparity between the two circuit types can be arbitrarily magnified by selecting KP problems with larger capacity bounds for the weights $w_i$. For instance, calculating the TTS difference at $n=15$ with a moderately realistic item weight distribution $w_i \sim U(1, 1000)$, the median TTS difference would ascend to a minimum of two orders of magnitude. Nevertheless, even for $n=15$ with $w_i \sim U(1,10)$, we are swiftly approaching the computational limits for simulating outcomes with QUBO-based circuits.

\section{Conclusion}
In this study, we analyze two approaches to solving combinatorial optimization problems with constraints within the framework of circuit-model quantum optimization based on adiabatic evolution. We argue that the prevalent traditional approach of recasting a constrained optimization problem into a quadratic Ising model using the QUBO formalism is highly inefficient; moreover, in addition to being costly in terms of the required computational resources, it often results in ill-behaved optimization landscapes. To address these challenges, we develop an alternative 
adiabatic quantum computation scheme that is based on the principle of Lagrangian duality from classical theory of optimization. First, we establish a theoretical connection between $\epsilon$-optimality, complementary slackness and the DAQC framework. Specifically, we show that the probability of sampling the optimal solution and the total evolution time are related by $\epsilon$. Secondly, by employing quantum circuits derived from a Trotterized approximation of adiabatic evolution, we implement a discrete version of adiabatic quantum computation (DAQC). We show that a DAQC scheme based on Lagrangian dual can be much more efficient than a DAQC scheme based on the QUBO approach. In particular, the Lagrangian dual approach results in significantly lower circuit connectivity and circuit depth compared to the latter. Furthermore, it also eliminates the need for logarithmic overheads in qubit count due to constraint bounds. We have conducted a numerical benchmarking study, which suggests that the proposed Lagrangian dual DAQC circuits yield a quadratic improvement over QUBO-based circuits in terms of adiabatic evolution time and circuit depth while achieving a comparable performance in solving a given constrained optimization problem. We further found that while LD-DAQC's performance is consistent regardless of a problem's coefficient values, the QUBO-based approach manifests a polylogarithmic growth in $R_{99}$ as the coefficient value range expands.

While the present study relies on discretized adiabatic scheduling, 
our proposed methodology is also applicable to variational algorithms  
and diabatic state preparation schemes, such as QAOA. In particular, 
the adiabatic schedule can be used as an initial guess and further refined through parameter optimization using a classical optimizer in a feedback loop. Such an approach could potentially further decrease the circuit depth or increase the algorithmic performance in terms of the expected time required to find an optimal solution for a given problem with high confidence.

Future investigations may address the potential for further performance or resource requirements improvements through a clever choice of the coupling terms $K_{i,j}X_iX_j$ in the mixing Hamiltonian given in \cref{eq:general_hamiltonian}. For example, efficiently computable metrics for linear problems, such as $\epsilon$, could be integrated into the coefficients $K_{i,j}$. Alternatively, leveraging knowledge of an $\epsilon$-optimal solution, which is readily calculable for these problems, might inform the design of an initial Hamiltonian or an initial state. 

Finally, we anticipate that the Lagrangian dual adiabatic quantum computation methodology proposed in the present work could also be applied to the more general framework of {\em mixed-integer programming} (MIP)~\cite{floudas1995nonlinear,li2006nonlinear,burer2012non}. Indeed, a recent study proposed an approach to MIP based on continuous-variable adiabatic quantum computation (CV-AQC) on a set of bosonic quantum field modes (qumodes)~\cite{khosravi2021mixed}. In that study, inequality constraints were handled by including quadratic penalty terms in the problem Hamiltonian along with additional auxiliary qumodes to encode the slack variables. It would be interesting to investigate if the Lagrangian dual methodology could potentially improve CV-AQC schemes in solving MIP problems with constraints.

\section{Acknowledgements}
The authors thank Pooya Ronagh for his critical feedback and suggestions that  greatly improved this work.
EG would also like to express his gratitude to Ben Adcock for his invaluable guidance, constant support, patience, and immense knowledge throughout the course of this research. EG also gratefully acknowledges the support from the Mitacs Accelerate Fellowship program.

\appendix
\section{\label{ap:duality_gap_overlap}Adiabatic Evolution via Lagrangian Duality Gap Analysis}
This section presents a detailed analysis showing how a small duality gap and the adiabatic transition performed in finite time can allow for sampling the solution of the primal problem with high probability. We develop the primal and dual optimization problems and incorporate the duality gap information into the adiabatic evolution.

Let $f(x):\{0,1\}^n \rightarrow \mathbb{R}$ and $g_i(x):\{0,1\}^n \rightarrow \mathbb{R}$ for $i=1,\ldots, m$. Consider the following optimization problem:
\begin{align}
    P = &\min_x f(x) \nonumber \\
    &\text{subject to} \nonumber\\
    &g_i(x) \geq c_i \text{ for } i =1, \ldots , m \nonumber \\
    &x \in \{0,1\}^n.
    \label{eq-app:primal-problem}
\end{align}
Then, the Lagrangian dual problem is given by:
\begin{equation}\label{eq:appendix_dual}
    D = \max_{\lambda \in \mathbb{R}^m_{+}} \min_x \left ( f(x) + \sum_{i=1}^m \lambda_i (c_i - g_i(x)) \right).
\end{equation}
Let $L(x)$ denote the Lagrangian function
\begin{equation}
    L(x) = f(x) + \sum_{i=1}^m \lambda_i^* (c_i - g_i(x)),
\end{equation}
where $\lambda_i^*>0$ are fixed and optimal for \cref{eq:appendix_dual}.

Let $x_P$ be the optimal solution of the primal problem, and let $x_D$ be the optimal solution of $L(x)$. Since $x_D$ is the optimal solution of the dual problem and $x_P$ is the optimal solution of the primal problem, we have:
\begin{equation}
    L(x_D) \leq L(x_P) \text{ and } f(x_P) \leq f(x_D).
\end{equation}
In the physical context, the problem Hamiltonian associated with the Lagrangian function $L(x)$ has the eigenstate $\ket{x_P}$ at a higher energetic level, which might be easier to reach when the adiabatic transition time is finite.

It is straightforward to show that
\begin{equation}\label{eq:eps_gap}
    L(x_P) - L(x_D) \leq \epsilon,
\end{equation}
where 
\begin{equation}
    \epsilon := \sum_{i=1}^m \lambda_i^* (g_i(x_D) - c_i).
\end{equation}
To show this, we have:
\begin{align*}
    L(x_P) - L(x_D) &= f(x_P) + \sum_{i=1}^m \lambda_i (c_i - g_i(x_P)) \\
    & \ \ \ \ - f(x_D) - \sum_{i=1}^m \lambda_i (c_i - g_i(x_D))\\
    &\leq f(x_P) - f(x_D) \\& \ \ \ \ - \sum_{i=1}^m \lambda_i (c_i - g_i(x_D)) \\
    &= f(x_P) - f(x_D) + \epsilon \\
    &= \epsilon - (f(x_D) - f(x_P))\\
    &\leq \epsilon,
\end{align*}
where we used the fact that $x_P$ must satisfy the constraints $g_i(x_P) \geq c_i$ according to \cref{eq-app:primal-problem}.
Thus, the energy gap between the eigenvalues $L(x_P)$ and $L(x_D)$ of the problem Hamiltonian representing $L(x)$ is bounded by $\epsilon$. If the evolution is not too slow and $\epsilon$ is sufficiently small, there can be a significant probability for the system to transition to a  higher energy level that is $\epsilon$-close. 

Next, we show that the probability amplitude of an energetically $\epsilon$-close state can be made arbitrarily close to 1 by an appropriate choice of evolution time $T$ and an initial Hamiltonian $H_{\text{init}}$. That is, even though the ground state $\ket{x_D}$ of the problem Hamiltonian is not the solution to the primal problem, the state $\ket{\psi(T)}$ prepared by LD-DAQC can have a significant overlap with both $\ket{x_D}$ and the solution of the primal problem $\ket{x_P}$.

Consider the adiabatic evolution of a quantum state $\ket{\psi(t)}$ such that $\ket{\psi(0)}$ is an eigenstate of the initial Hamiltonian $H_{\text{init}}$. We write $\ket{\psi(t)}$ in the instantaneous eigenbasis of the total Hamiltonian $H(t)$:
\begin{equation}
    \ket{\psi(t)} = \sum_k c_k(t)\ket{\psi_k(t)},
\end{equation}
where the instantaneous eigenstates $\ket{\psi_k(t)}$ are labeled such that the eigenvalues are ordered as $E_0(t) < \cdots < E_k(t) \ldots < E_{2^n-1}(t)$. Furthermore, we exploit the freedom in choosing the global phase of the eigenstate, $\ket{\psi_k(t)}$, such that $\langle \psi_k(t) | \dot{\psi}_k(t) \rangle = 0$. This is done to simplify subsequent computations.

We now create a system of coupled differential equations that describe the evolution of the coefficients $c_k(t)$. To this end, differentiate $\ket{\psi(t)}$:\
\begin{equation}\label{eq:deriv_psi}
    \ket{\dot{\psi}(t)} = \sum_k \left( \dot{c}_k(t)\ket{\psi_k(t)} + c_k(t)\ket{\dot \psi_k(t)} \right).
\end{equation}
Substitute this into the Schr\"odinger equation:
\begin{align}
    \mathrm{i} \ket{\dot{\psi}(t)} = H(t) \ket{\psi(t)}.
\end{align}
Expand the above equation using \cref{eq:deriv_psi}:
\begin{align}
    \mathrm{i} \sum_k \left( \dot{c}_k(t)\ket{\psi_k(t)} + c_k(t)\ket{\dot \psi_k(t)} \right) = \nonumber \\
    \sum_k c_k(t)E_k(t)\ket{\psi_k(t)}.
\end{align}
To isolate a coefficient $c_j(t)$, we project both sides of the equation onto eigenstate $\ket{\psi_j(t)}$ and obtain:
\begin{align}\label{eq:full_system_coeffs}
    \textrm{i} \dot c_j(t) = c_j(t)E_j(t) - \textrm{i} &\sum_k c_k(t) \braket{\psi_j(t)}{\dot \psi_k(t)} \nonumber \\
    &\text{ for } j=0,\ldots, 2^n-1.
\end{align}
We have obtained a system of coupled differential equations that describe the evolution of coefficients $c_j(t)$. We will refer to the term $\braket{\psi_j(t)}{\dot \psi_k(t)}$ as \textit{non-adiabatic couplers}. When the evolution is fully adiabatic, these couplers are zero, and the system of differential equations is uncoupled:
\begin{equation}\label{eq:full_system_coeffs_adiabat}
    \textrm{i} \dot c_j(t) = c_j(t)E_j(t) \text{ for } j=0,\ldots, 2^n-1.
\end{equation}
Given the initial condition $c_0(0)=1$ and $c_j(0)=0$ for $j>0$, it is easy to deduce that $|c_0(t)|^2=1$ and $|c_j(t)|^2 = 0$  for all $t$ and $j >0$. That is, whenever the system evolves adiabatically, the state $\ket{\psi(t)}$ is $c_0(t)\ket{\psi_0(t)}$. In other words, the system stays in its lowest instantaneous energy eigenstate. This is because the differential equation $\mathrm{i} \dot c_0(t) = c_0(t)E_0(t)$, subject to the initial condition $c_0(0)=1$, has the following solution:
\begin{equation}
    c_0(t) = e^{-\rm{i} \int_0^t E_0(s)ds}.
\end{equation}
However, if the evolution is not entirely adiabatic, the non-adiabatic couplers $\braket{\psi_j(t)}{\dot \psi_k(t)}$ are non-zero. Hence, the coefficients $c_j(t)$ in \cref{eq:full_system_coeffs} have highly nontrivial dynamics.

Generally, the system of differential equations in \cref{eq:full_system_coeffs} is intractable and an analytical analysis is not possible. However, it is possible to greatly simplify the analysis by making some reasonable assumptions.

First, we rewrite the \cref{eq:full_system_coeffs} in terms of $\dot H(t)$ and energy gaps. It can be shown that, for $j \neq k$, a non-adiabatic coupler can be expressed as:
\begin{align}
    \braket{\psi_j(t)}{\dot \psi_k(t)} &= -\frac{\bra{\psi_j(t)}\dot H(t)\ket{\psi_k(t)}}{E_j(t)-E_k(t)}.
\end{align}
Let us compute the explicit form of $\dot H(t)$. For the simplicity of the analysis, we choose a linear adiabatic schedule. Then the total Hamiltonian is given by
\begin{equation}
    H(t) = \left(1 - \frac{t}{T}\right)H_{\text{init}} + \frac{t}{T} H_P .
\end{equation}
Its derivative is:
\begin{equation}
    \dot H(t) = -\frac{1}{T} H_{\text{init}} + \frac{1}{T}H_P = \frac{1}{T}(H_P - H_{\text{init}}). 
\end{equation}
Hence,
\begin{equation}
    \bra{\psi_j(t)} \dot H(t) \ket{\psi_k(t)} = \bra{\psi_j(t)} \frac{1}{T}(H_{\text{P}} - H_{\text{init}}) \ket{\psi_k(t)},
\end{equation}
and the coupler becomes
\begin{equation}
    \braket{\psi_j(t)}{\dot \psi_k(t)} = -\frac{\bra{\psi_j(t)}(H_{\text{P}} - H_{\text{init}}) \ket{\psi_k(t)}}{T(E_j(t)-E_k(t))}.
\end{equation}
We define $D_{jk}(t)$ to be the numerator of the coupler:
\begin{equation}
    D_{jk}(t) := \bra{\psi_j(t)}(H_{\text{P}} - H_{\text{init}}) \ket{\psi_k(t)}.
\end{equation}
If the adiabatic transition time $T$ tends to infinity, the coupler tends to zero and we obtain the adiabatic evolution described by \cref{eq:full_system_coeffs_adiabat}. However, if the energy gap between $E_j(t)$ and $E_k(t)$ is small, then the coupler cannot be ignored. Given the new formulation of the coupler, the system of differential equations is given by:
\begin{align}
    \label{eq:full_system_coeffs_gap}
    \textrm{i} \dot c_j(t) = c_j(t)E_j(t) + \textrm{i} \sum_{k \neq j} c_k(t) \frac{D_{jk}(t) }{T\left(E_j(t)-E_k(t)\right)}\\
    \text{ for } j=0,\ldots,2^n-1 \nonumber.
\end{align}

We now focus on the coefficient $c_0(t)$ of the instantaneous ground state $\ket{\psi_0(t)}$. In the case of an ideal adiabatic evolution, $\ket{\psi_0(t)}$ would evolve into $\ket{x_D}$. The differential equation of $c_0(t)$ is given by:
\begin{equation}
    \label{eq:zero_system_coeffs}
    \textrm{i} \dot c_0(t) = c_0(t)E_0(t) - \textrm{i}\sum_{k\neq0} c_k(t) \frac{D_{0k}(t)}{T\left(E_k(t)-E_0(t)\right)}.
\end{equation}
Our first simplifying assumption is to drop the coupling terms with a large energy gap $E_k(t)-E_0(t)$. This is reasonable, as a large gap along with suitably chosen $T$ (for the evolution to be adiabatic) significantly diminish the contribution of the term $D_{0k}(t)/ T(E_k(t)-E_0(t))$. Specifically, we will ignore terms that exceed a certain threshold. To this end, we fix some time $t^*$ and work with the following index set:
\begin{equation}\label{eq:filtering_ineq}
    Z := \left \{ k \in \mathbb{N} \ | \ E_k(t)- E_0(t) \leq \epsilon \text{ for } t^* < t < T \right \}.
\end{equation}
Let $z := \max Z$, then for all $k \leq z$ eigenvalues $E_k(t)$ satisfy the inequality $E_k(t)- E_0(t) \leq \epsilon$. If $\epsilon$ is sufficiently small, then we expect $z$ to be small as well. Hence, the states $\ket{\psi_k(t)}$ for $k \in Z$ and $t > t^*$ are energetically close to the instantaneous ground state $\ket{\psi_0(t)}$, and therefore, their coupling terms $D_{0k}(t)/T\left(E_k(t)-E_0(t)\right)$ for $k=1,\ldots, z$ have the most significant contribution. Thus, after dropping the terms with $k \notin Z$, we consider the following evolution:
\begin{equation}\label{eq:reduced_system}
    \textrm{i} \dot c_0(t) = c_0(t)E_0(t) - \textrm{i}\sum_{0 < k \leq z} c_k(t) \frac{D_{0k}(t) }{T\left(E_k(t)-E_0(t)\right)}.
\end{equation}

This simplification is still not sufficient to make the analysis tractable. The second simplification is achieved by a specific choice of the coefficients of the Hamiltonian $H_{\text{init}}$. First, note that at time $T$, the eigenstates $\ket{\psi_k(T)}$ are the computational basis states, as $H_{\text{P}}$ is diagonal in the computation basis. Therefore, for $k \neq 0$, we have $D_{0k}(T) = \bra{\psi_0(T)}(H_{\text{P}}-H_{\text{init}})\ket{\psi_k(T)} =-\bra{\psi_0(T)}H_{\text{init}}\ket{\psi_k(T)} \equiv -[H_{\text{init}}]_{0k}$. Because we are free to choose any initial Hamiltonian as long as it does not commute with $H_{\text{P}}$, we suppose that $H_{\text{init}}$ was chosen such that $[H_{\text{init}}]_{0k}=0$ for $k=1,\ldots, z-1$ and $[H_{\text{init}}]_{0z} \neq 0$. It follows that, for $k=1,\ldots, z-1$, $D_{0k}(t)$ vanishes as $t$ approaches $T$, i.e. $D_{0k}(T) = [H_{\text{init}}]_{0k} = 0$. Therefore, we will drop them to further simplify the computation. The only coefficent that does not vanish for all $t$ is $D_{0z}(t)$. The values of $D_{0z}(t)$ are not known and intractable for all $t < T$, however, at $t=T$, we have $D_{0z}(T) = -[H_{\text{init}}]_{0z}$. Therefore, we will approximate $D_{0z}(t)$ with $-[H_{\text{init}}]_{0z}$. Then the differential equation simplifies to:
\begin{equation}
    \textrm{i}  \dot{\tilde{c}}_0(t) = \tilde c_0(t)E_0(t) + \textrm{i} \tilde c_z(t) \frac{[H_{\text{init}}]_{0z} }{T\left(E_z(t)-E_0(t)\right)},
\end{equation}
where $\tilde{c}_0(t)$ and $\tilde{c}_z(t)$ are approximations of $c_0(t)$ and $c_z(t)$, respectively. We apply the analogous simplifying assumptions to the differential equation describing the evolution of the coefficient $c_z(t)$. Specifically, we only consider the coupling terms with indices $k < z$, and set $[H_{\text{init}}]_{zk}=0$ for $0 < k < z$. Then, we obtain the following coupled differential equations:
\begin{align}\label{eq:final_system}
    \dot{\tilde{c}}_0(t) &= -\textrm{i} \tilde c_0(t)E_0(t) + \tilde c_z(t) \frac{A}{T\Delta E(t)}, \nonumber\\
    \dot{\tilde c}_z(t) &= -\textrm{i} \tilde c_z(t)E_z(t) - \tilde c_0(t) \frac{A}{T\Delta E(t)}.
\end{align}
In the above, we define $A:= [H_{\text{init}}]_{0z} = [H_{\text{init}}]_{z0}$, and $\Delta E(t) := E_z(t)-E_0(t)$. From computation in \cref{eq:eps_gap} we know that $\Delta E(T) = L(x_P)-L(x_D) \leq \epsilon$. In other words, by the end of the evolution, the gap is bounded by $\epsilon$. Because we do not have knowledge of \mbox{$\Delta E(t) = E_z(t)-E_0(t)$}, and given that the transition schedule is linear, we interpolate $\Delta E(t)$ with a linear function:
\begin{equation}
    p(t) = r(T - t) + \epsilon.
\end{equation}
In the above, $r$ is some constant such that $p(0) = \Delta E(0)$, and $p(T) = \epsilon \geq \Delta E(T)$. Then, the system of differential equations is:
\begin{align}\label{eq:final_system_with_approximated_gap}
    \dot{\tilde{c}}_0(t) &= -\textrm{i} \tilde c_0(t)E_0(t) + \tilde c_z(t) \frac{A}{T p(t)}, \nonumber\\
    \dot{\tilde{c}}_z(t) &= -\textrm{i} \tilde c_z(t)E_z(t) - \tilde c_0(t) \frac{A}{T p(t)}.
\end{align}
We note that $A$ defined as $[H_{\text{init}}]_{0z} \neq 0$ depends on the initial Hamiltonian which can be freely chosen. Therefore, we can choose $A$ and $T$ such that the second terms on the right hand side of \cref{eq:final_system_with_approximated_gap} dominate the dynamics. Hence, we can further simplify the system to:
\begin{align}\label{eq:final_system_simplified}
    \dot{\tilde c}_0(t) &=\tilde c_z(t) \frac{A}{T p(t)}\nonumber\\
    \dot{\tilde c}_z(t) &= -\tilde c_0(t) \frac{A}{T p(t)}.
\end{align}
We can rewrite this in a matrix form:
\begin{equation}
    \begin{pmatrix}
        \dot{\tilde c}_0(t) \\
        \dot{\tilde c}_z(t) \\
    \end{pmatrix}=
    \begin{pmatrix}
        0 &&  \frac{A}{T p(t)}\\
        - \frac{A}{T p(t)} && 0
    \end{pmatrix}
    \begin{pmatrix}
        \tilde{c}_0(t) \\
        \tilde{c}_z(t)
    \end{pmatrix}.
\end{equation}
Given the initial conditions $c_0(0) = 1$ and $c_z(0) = 0$, we obtain the following solutions:
\begin{align}\label{eq:approximate_probabilities}
    \tilde{c}_0(t) &= \cos g(t),\nonumber \\
    \tilde{c}_z(t) &= - \sin g(t),
\end{align}
where $g(t)$ is defined as:
\begin{align}
    g(t) := \frac{A \ln(Tp(0))- \ln(Tp(t))}{Tr}.
\end{align}
Specifically, a time $T$ we have
\begin{equation}
    g(T) = \frac{A}{Tr}\ln(1 + \frac{Tr}{\epsilon}).
\end{equation}
Hence, the probability amplitudes at the final time $T$ are:
\begin{align}
    \tilde{c}_0(T) &= \cos \left(\frac{A}{Tr}\ln(1 + \frac{Tr}{\epsilon})\right), \nonumber \\
    \tilde{c}_z(T) &= - \sin \left(\frac{A}{Tr}\ln(1 + \frac{Tr}{\epsilon})\right).
\end{align}
An appropriate choice of $T, r$ and $A$ makes it possible to achieve $|\tilde{c}_z(T)|^2=1$.

\begin{figure}[t]
\centering
\includegraphics[width=0.45\textwidth]{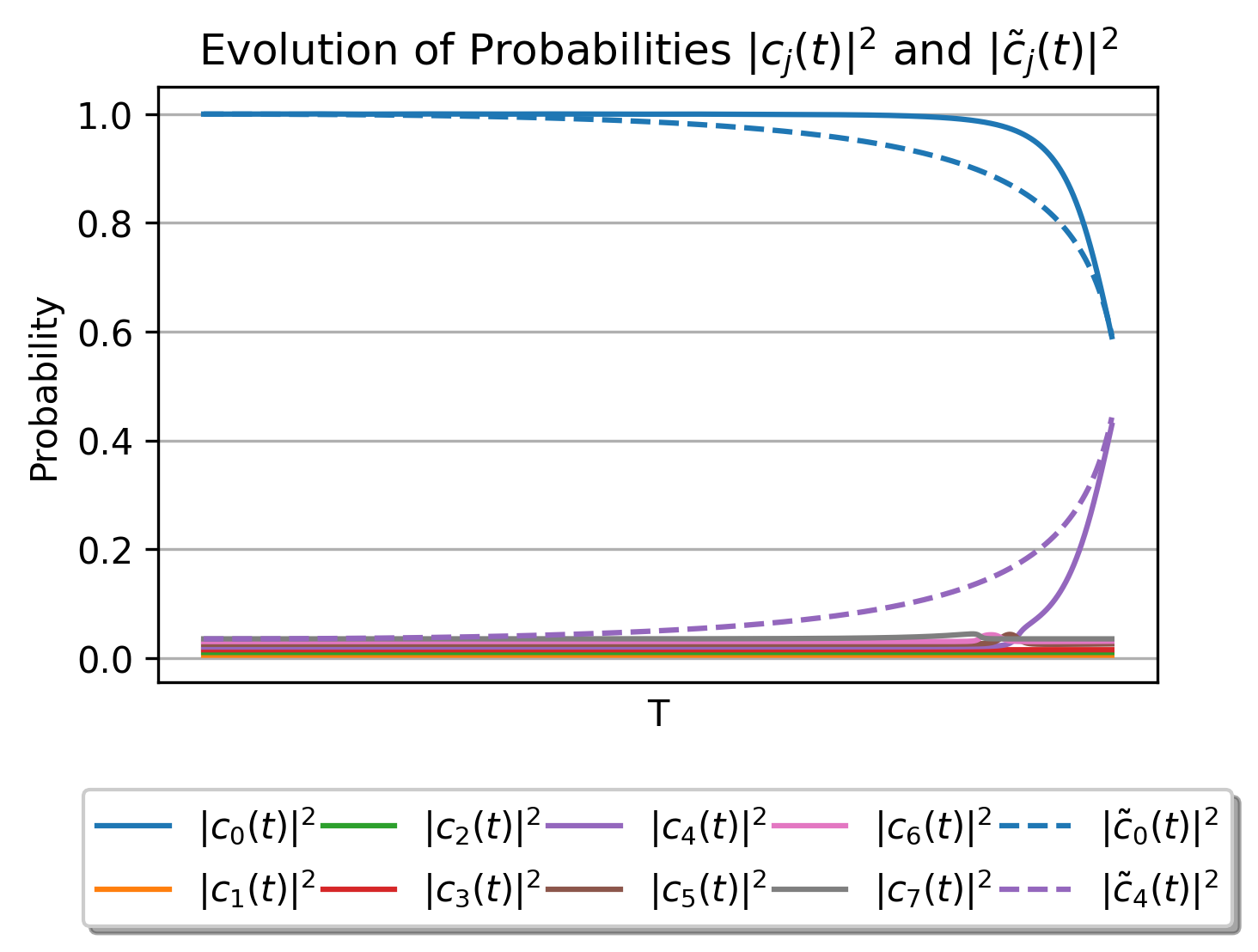}
\caption{
    The evolution of the instantaneous probabilities $|c_j(t)|^2$ for $j=0,\ldots, 7$ (solid lines) is depicted under the system of differential equations presented in \cref{eq:full_system_coeffs_gap}, alongside the evolution of the approximations (dashed lines) $|\tilde{c}_0(t)|^2$ and $|\tilde{c}_z(t)|^2$, where $z=4$. It is observed that $|c_z(t)|^2$, which represents the probability of sampling $\ket{\psi_z(t)}$, increases over time. By the end of the evolution, we note that $|c_z(T)|^2 \approx |\tilde{c}_z(T)|^2 > 0.4$. Consequently, $|\tilde{c}_z(T)|$ closely approximates the probability of sampling the solution to the primal problem $\ket{\psi_z(T)}= \ket{x_P}$.
    }\label{fig:evo_of_coeffs_c(t)}
\end{figure}

\textbf{Numerical investigation} Numerical experiments have shown that, despite the use of simplifying approximations, the dynamics of the probabilities $|\tilde{c}_0(t)|^2$ and $|\tilde{c}_z(t)|^2$ in \cref{eq:approximate_probabilities} closely follow that of $|c_0(t)|^2$ and $|c_z(t)|^2$ given by the system of differential equation in \cref{eq:full_system_coeffs_gap} where we have not made any simplifications. Specifically, for any given problem Hamiltonian $H_{\text{P}}$ and an initial Hamiltonian $H_{\text{init}}$ constructed as previously described, it is feasible to achieve $|\tilde{c}_0(T)|^2 \approx |c_0(T)|^2$ and $|\tilde{c}_z(T)|^2 \approx |c_z(T)|^2$ by adjusting the parameter $r$, which determines the rate of change of $p(t)$. For example, \cref{fig:evo_of_coeffs_c(t)} illustrates the evolution of the instantaneous probabilities $|c_{j}(t)|^2$ for $j=0,\ldots, 7$, as determined by the system of differential equations in \cref{eq:full_system_coeffs_gap}. In this specific case, we have $\epsilon = 2.5$ and $z = 4$. It is noteworthy that $|\tilde{c}_z(T)|^2 \approx |c_z(T)|^2$ and $|\tilde{c}_0(T)|^2 \approx |c_0(T)|^2$. Additionally, it is observed that at the beginning of the evolution, the probability of sampling the instantaneous ground state $\ket{\psi_0(t)}$ is almost one, as expected. However, over time, due to the choice of $T$ and the construction of the initial Hamiltonian, the probability of sampling the state $\ket{\psi_z(t)}$ increases, while the probabilities of other states either remain constant or decrease. By the time $T$, when $\ket{\psi_z(T)} = \ket{x_P}$, the probabilities satisfy $|c_z(T)|^2 \approx |\tilde{c}_z(T)|^2 > 0.4$. Thus, by the end of the evolution, the approximation $|\tilde{c}_z(T)|^2$ successfully models the probability of sampling the solution to the primal problem $\ket{\psi_z(T)}=\ket{x_P}$.

\section{\label{ap:complexity}Analysis of Runtime Complexity for LD-DAQC and QUBO-DAQC Circuits}
Recall that one-qubit gates acting on different qubits can be executed in parallel, whereas any pair of two-qubit gates sharing at least one qubit can only be executed sequentially. As a result, highly parallelizable circuits have shorter runtime and are less prone to accumulation of errors.
This section analyzes the runtime complexity of the Lagrangian- and QUBO-based circuits for a KP. In Section~\ref{sec:metrics} we stated that the LD-DAQC circuit's runtime $t_{\text{ss}}$ is $O(p)$ where $p$ is the number of layers. In contrast, the QUBO-DAQC circuit has runtime complexity $O\left((n + \log_2(c)) p\right)$ where $n$ is the number of variables in a KP and $c$ is the capacity.

\begin{figure}[h!]
\centering
\includegraphics[width=0.47\textwidth]{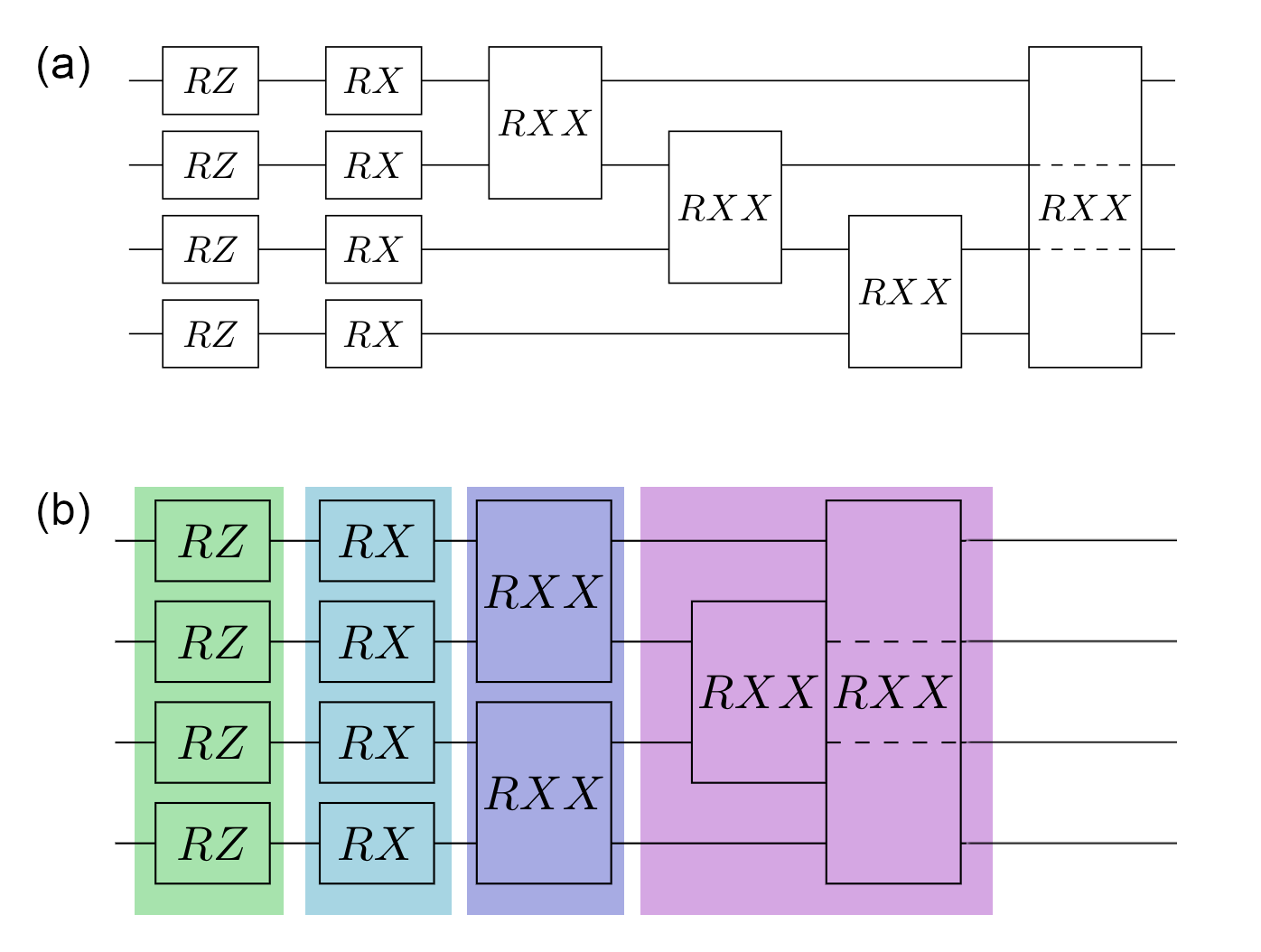}
\caption[]{(a) The $k$th layer of the LD-DAQC circuit. By the choice of the coefficients $K_{i,j}$ in \cref{eq:general_hamiltonian} we get the initial Hamiltonian \cref{eq:mixer}, which in turn yields $RXX$ gates with a closed chain-like connectivity.
(b) Equivalent rearrangement of the $k$th layer. Gates are partitioned into four sublayers (coloured rectangles). Gates belonging to the same sublayer do not share qubits and can be applied in parallel in a single step. In total it takes four consecutive executions to apply the $k$th layer for any even $n$.}\label{fig:lagrangian_circuit_rearrange}
\end{figure}

Let us first look at the LD-DAQC circuit given in \cref{eq:lag_approx_U}-\cref{eq:gates_problem}. To show that the runtime $t_{\text{ss}}$ is $O(p)$ it is sufficient to examine any single layer of the circuit illustrated in Fig~\ref{fig:daqc_circuit}. The $k$th layer has the following form:
\begin{align}\label{eq:kth_layer}
    &\prod_{j=1}^n \exp \left \{ \mi \gamma_k X_j X_{j+1} \right \} \times \prod_{j=1}^n \exp \left \{ \mi \gamma_k X_j \right \} \nonumber\\
    & \times \prod_{j=1}^n \exp \left \{ -\mi \beta_k \left( v_j - \lambda(\Delta t k)w_j \right)Z_j \right \}.
\end{align}
We claim that the runtime of a single layer $k$ is constant. Intuitively, we can rearrange the gates of the $k$th layer into at most five sublayers such that gates belonging to a sublayer can be applied in parallel. Each sublayer has the runtime $O(1)$. Fig~\ref{fig:lagrangian_circuit_rearrange} illustrates a 4-qubit $k$th layer of a LD-DAQC circuit and its equivalent rearrangement into four sublayers that can be applied consecutively. For all $j = 1, \ldots, n$ the unitary matrix $\exp \left \{ \mi \gamma_k X_j \right \}$ represents a one-qubit $RX_j$ gate acting on the qubit $j$. Therefore, their product can be applied in parallel in one step. Similarly, the matrix $\exp \left \{ -\mi \beta_k \left( v_j - \lambda(\Delta t k)w_j \right)Z_j \right \}$ represents a one-qubit $RZ_j$ gate acting on the qubit $j$. Hence, the product of all $RZ$ gates can be applied in parallel in one step. We conclude that the number of steps required to apply all $RX$ and $RZ$ gates is independent of $n$, i.e., it is constant per layer $k$.

We now examine the two-qubit gates $RXX_{j,j+1}$ given by $\exp \left \{ \mi \gamma_k X_j X_{j+1} \right \}$. We note that their product forms a closed chain where the qubit $j$ is connected to the qubit $j+1$, recalling that $n + 1 := 1$ (periodic boundary condition). Since for any $i$ and $j$ the matrices $\exp \left \{ \mi \gamma_k X_j X_{j+1} \right \}$ and $\exp \left \{ \mi \gamma_k X_i X_{i+1} \right \}$ commute we can rearrange their order into several sublayers so that each layer could be applied in parallel. Suppose we have an even number of qubits, i.e., $n = 2m$ for some $m \in \mathbb{N}$. Then the first sublayer will contain $m$ gates that do not share any qubits and the second sublayer will contain the other half of $m$ gates that do not share any qubits. Hence, the two sublayers of $RXX$ gates can be executed in 2 steps. If $n$ is odd, then we get a third sublayer containing a single $RXX$ gate. Hence, the three sublayers of $RXX$ gates can be executed in 3 steps. It follows that regardless of the value of $n$ it takes at most 3 steps to apply all the $RXX$ gates. Therefore, in total it takes at most 5 steps to apply the $RX, RZ$ and $RXX$ gates. Hence each DAQC layer $k$ in \cref{eq:kth_layer} can be applied in constant time. Since there are $p$ layers the runtime $t_{\text{ss}}$ is $O(p)$.

For the QUBO circuit, on the other hand, the runtime depends on the capacity $c$, the number of KP variables $n$, and the number of layers $p$. Due to the quadratic penalty in the objective function \cref{eq:qubo} the $k$th DAQC layer has $\binom{n + \lfloor \log_2(c) \rfloor}{2}$ two-qubit gates $RZZ$. Since $RZZ$ gates commute, the best arrangement of the $RZZ$ gates yields at least $n + \lfloor \log_2(c) \rfloor - 1$ sublayers that can be executed consecutively. To see this, we note that each sublayer can accommodate at most $(n + \lfloor \log_2(c) \rfloor) / 2$ two-qubit gates such that no qubit is shared between all the gates in the sublayer. Therefore it takes at least $n + \lfloor \log_2(c) \rfloor - 1$ sublayers to place all the $RZZ$ gates. This gives $(n + \lfloor \log_2(c) \rfloor - 1)(n + \lfloor \log_2(c) \rfloor) / 2 = \binom{n + \lfloor \log_2(c) \rfloor}{2}$. Since the sublayers containing $RZ$ and $RX$ can be applied in constant time, the runtime of a $k$th layer is of order $O(n + \lfloor \log_2(c) \rfloor)$. Since the QUBO circuit has $p$ layers, its runtime is $t_{\text{ss}} = O(p(n + \lfloor \log_2(c) \rfloor))$.

\bibliography{main.bib}
\end{document}